\documentclass[12pt]{amsart}
\usepackage{amssymb}
\usepackage{amsmath}
\usepackage[dvips]{color}
\usepackage{graphicx}
\usepackage{epsfig}
\usepackage[dvips,line,arrow,curve,rotate]{xy}
\usepackage{color}
\newtheorem{theorem}{Theorem}[section]
\newtheorem{lemma}[theorem]{Lemma}
\newtheorem{proposition}[theorem]{Proposition}
\newtheorem{corollary}[theorem]{Corollary}

\theoremstyle{definition}
\newtheorem{example}[theorem]{Example}

\newtheorem{definition}[theorem]{Definition}

\theoremstyle{remark}

\newtheorem{remark}[theorem]{Remark}

\DeclareMathOperator{\Dl}{\mathcal{L}({\it D})}
\DeclareMathOperator{\Dm}{\mathcal{M}({\it D})}
\DeclareMathOperator{\Dr}{\mathcal{R}({\it D})}

\newlength{\cellsz}
\newcounter{cellsize}
\newcommand{\setcellsize}[1]{%
  \setcounter{cellsize}{#1}%
  \setlength{\cellsz}{\value{cellsize}\unitlength}}%
\setcellsize{13}%
\newcommand\cellify[1]{\def\thearg{#1}\def\nothing{}%
\ifx\thearg\nothing \vrule width0pt height\cellsz depth0pt\else
\hbox to 0pt{{\begin{picture}(\value{cellsize},\value{cellsize})
  \put(0,0){\line(1,0){\value{cellsize}}}
  \put(0,0){\line(0,1){\value{cellsize}}}
  \put(\value{cellsize},0){\line(0,1){\value{cellsize}}}
  \put(0,\value{cellsize}){\line(1,0){\value{cellsize}}} \end{picture} \hss}}\fi%
\vbox to \cellsz{ \vss \hbox to \cellsz{\hss$#1$\hss} \vss}}
\newcommand\tableau[1]{\vcenter{\vbox{\let\\\cr
\baselineskip -16000pt \lineskiplimit 16000pt \lineskip 0pt
\ialign{&\cellify{##}\cr#1\crcr}}}}
\newcommand\tabl[1]{\vtop{\let\\\cr
\baselineskip -16000pt \lineskiplimit 16000pt \lineskip 0pt
\ialign{&\cellify{##}\cr#1\crcr}}}

\def\noXv{\,\,\lower2pt\hbox{
\begin{picture}(0,0)%
\includegraphics{figNoX.pstex}%
\end{picture}%
\setlength{\unitlength}{1973sp}%
\begingroup\makeatletter\ifx\SetFigFont\undefined%
\gdef\SetFigFont#1#2#3#4#5{%
  \reset@font\fontsize{#1}{#2pt}%
  \fontfamily{#3}\fontseries{#4}\fontshape{#5}%
  \selectfont}%
\fi\endgroup%
\begin{picture}(316,316)(293,-969)
\end{picture}
}}
\def\noXvsp{\,\,\lower2pt\hbox{\begin{picture}(0,0)%
\includegraphics{noXvsp.pstex}%
\end{picture}%
\setlength{\unitlength}{1657sp}%
\begingroup\makeatletter\ifx\SetFigFontNFSS\undefined%
\gdef\SetFigFontNFSS#1#2#3#4#5{%
  \reset@font\fontsize{#1}{#2pt}%
  \fontfamily{#3}\fontseries{#4}\fontshape{#5}%
  \selectfont}%
\fi\endgroup%
\begin{picture}(699,474)(7189,-3223)
\end{picture}%
}}

\begin{document}

\author[T. K. Petersen]{T. Kyle Petersen}

\author[P. Pylyavskyy]{Pavlo Pylyavskyy}

\author[B. Rhoades]{Brendon Rhoades}

\title{Promotion and cyclic sieving via webs}

\begin{abstract}
We show that Sch\"utzenberger's promotion on two and three row rectangular Young tableaux can be realized as cyclic rotation of certain planar graphs introduced by Kuperberg. Moreover, following work of the third author, we show that this action admits the cyclic sieving phenomenon. 
\end{abstract}

\maketitle

\section{Introduction}\label{sec:def}

Let us briefly recall some definitions; refer to \cite{St} for more details. A \emph{partition} $\lambda = (\lambda_1,\ldots,\lambda_k)$ of an integer $n$, written $\lambda \vdash n$, is multiset of positive integers whose sum is $n$, which by convention is written in weakly decreasing order. For every partition of $n$ we can draw an arrangement of $n$ boxes into left-justified rows of lengths $\lambda_1\geq \lambda_2\geq \cdots \geq\lambda_k$, called a \emph{Young diagram}. A \emph{semistandard Young tableau} is a way of filling the boxes in a Young diagram with positive integers so that the entries weakly increase in rows, strictly increase down columns. The \emph{type} of a semistandard Young tableau is the multiset of entries. A \emph{standard Young tableau} is a semistandard tableau of type $\{1,2,\ldots,n\}$, where $n$ is the number of boxes. Given a partition $\lambda$, let $SSYT(\lambda)$ denote the set of semistandard Young tableaux of shape $\lambda$, and similarly let $SYT(\lambda)$ denote the set of standard Young tableaux of shape $\lambda$. For example,
\[
\tableau{{1}&{1 }&{2 }&{4 } \\ {2 }&{2 }&{3 } \\ {3 }} \in SSYT((4,3,1)), \qquad \tableau{{1}&{2 }&{4 }&{8 } \\ {3 }&{6 }&{7 } \\ {5 }} \in SYT( (4,3,1) ).
\]
We denote the entry in row $a$, column $b$ of a tableau $T$, by $T_{a,b}$. Another way to describe a standard Young tableau is to write its \emph{Yamanouchi word}. The Yamanouchi word for a tableau $T \in SYT(\lambda)$, with $\lambda = (\lambda_1,\lambda_2, \ldots, \lambda_k) \vdash n$, is a word $w = w_1 \cdots w_n$ on the multiset \[\{1^{\lambda_1},2^{\lambda_2},\ldots, k^{\lambda_k}\} := \{ \underbrace{1,\ldots, 1}_{\lambda_1}, \underbrace{2,\ldots,2}_{\lambda_2}, \ldots, \underbrace{k,\ldots,k}_{\lambda_k}\},\] such that $w_i$ is the row in which $i$ is placed in $T$. For example,
\[
\tableau{{1}&{2 }&{4 }&{8 } \\ {3 }&{6 }&{7 } \\ {5 }} \leftrightarrow 11213221.
\] 
Notice that Yamanouchi words are characterized by the fact that in reading $w$ from left to right, there are never fewer letters $i$ than letters $(i+1)$. Given such a word $w$ we can associate a tableau $T(w)$ in a straightforward way. We say $w$ is \emph{balanced} if all distinct letters appear the same number of times. Balanced Yamanouchi words are in bijection with standard Young tableaux of rectangular shapes.

In this paper we will study the action of \emph{jeu-de-taquin promotion} on certain classes of tableaux.
Promotion was defined by Sch\"utzenberger as an action on posets \cite{Sch}, and has since appeared in a number of contexts, usually applied to tablueax, cf. \cite{H, Sh, Ste}. For our purposes, promotion is a bijection $p: SYT(\lambda) \to SYT(\lambda)$ defined as follows. 

\begin{definition}[Jeu-de-taquin promotion] 
Given a tableau $T$ in $SYT(\lambda)$ with $\lambda \vdash n$, form $p(T)$ with the following algorithm.
\begin{enumerate}
\item Remove the entry 1 in the upper left corner and decrease every other entry by 1. The empty box is initialized in position $(a,b) = (1,1)$.

\item Perform jeu de taquin:
\begin{enumerate}

\item If there is no box to the right of the empty box and no box below the empty box, then go to 3).

\item If there is a box to the right or below the empty box, then swap the empty box with the box containing the smaller entry, i.e., $p(T)_{a,b} := \min\{T_{a,b+1}-1, T_{a+1,b}-1\}$. Set $(a,b) := (a',b')$, where $(a',b')$ are the coordinates of box swapped, and go to 2a).

\end{enumerate}

\item Fill the empty box with $n$.

\end{enumerate}
\end{definition}

For example,
\[T = \tableau{{1}&{2 }&{4 }&{8 } \\ {3 }&{6 }&{7 } \\ {5 }} \quad \mapsto  \quad \tableau{{1}&{3 }&{6 }&{7 } \\ {2 }&{5 }&{8 } \\ {4 }} = p(T).
\]

\begin{remark}
We take care to point out that promotion should not be confused with the similarly defined, and more widely studied, action called \emph{evacuation} (or the ``Sch\"utzenberger involution," or, more confusingly, ``evacuation and promotion"), also defined in \cite{Sch}.
\end{remark}

As a permutation, promotion naturally splits $SYT(\lambda)$ into disjoint orbits. For a general shape $\lambda$ there seems to be no obvious pattern to the sizes of the orbits. However, for certain shapes, notably Haiman's ``generalized staircases" more can be said \cite{H} (see also Edelman and Greene \cite[Cor. 7.23]{EG}). In particular, rectangles fall into this category, with the following result.

\begin{theorem}[\cite{H}, Theorem 4.4] \label{thm:ordern}
 \item If $\lambda = (n,\ldots,n) \vdash N = bn$ is a rectangle, then $p^N(T) = T$ for all $T \in SYT(\lambda)$.
\end{theorem}

In this paper we will reinterpret the action of promotion on rectangular standard tableaux having two or three rows as a more elementary action on different sets of combinatorial objects.  
These alternative descriptions of the action of promotion will render Theorem 1.1 transparent for the cases $b = 2$ and $b = 3$.  In the case of $b = 2$ this interpretation was discovered by White \cite{W} and takes the form of a bijection from the set of standard tableaux of shape 2 by $n$ and the set of noncrossing matchings on $[2n]$ under which promotion on tableaux maps to rotation on matchings.
As we will show (Theorem \ref{thm:prom}), the case of $b = 3$ involves a bijection from the set of three row standard tableaux to a collection of combinatorial objects called \emph{$A_2$ webs} under which promotion maps to a combinatorial action called web rotation.

Let us now review the result for $b=2$ rows. Given a balanced Yamanouchi word $w=w_1\cdots w_{2n}$ on $\{1^n,2^n\}$, draw $2n$ vertices around the boundary of a disk, label them $1,\ldots,2n$ counterclockwise, and place $w_i$ at vertex $i$. Read the word and for every $2$ we encounter, draw a line between that vertex and the clockwise nearest $1$ that is not already matched with a $2$. To recover a Yamanouchi word from a noncrossing matching, traverse the disk counterclockwise, starting at the first vertex. On first encountering an edge, label the endpoint with a $1$, the second time label the endpoint with a $2$. Below are the five noncrossing matchings on six vertices labeled with the corresponding Yamanouchi words ($w_1$ is at 11 o'clock).

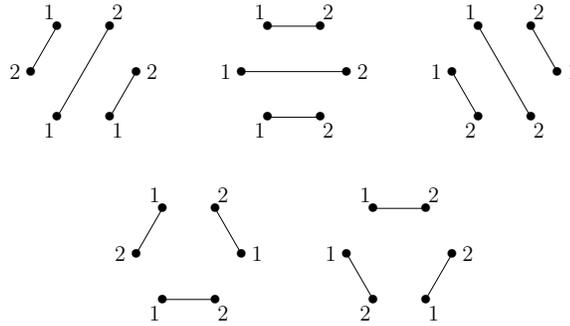
\[
\begin{xy}
<1cm,0cm>:
(0,0)*{\object[*.7]{
\xybox{
\begin{xy}
0;<1cm,0cm>:<.5cm,\halfrootthree cm>::
(2,-1)*{}, (0,7)*{},
(0,5); (0,6) **@{-}, (1,4); (1,6) **@{-}, (2,4); (2,5) **@{-},
(4,5); (6,5) **@{-}, (5,4); (6,4) **@{-}, (4,6); (5,6) **@{-},
(8,5); (9,4) **@{-}, (8,6); (10,4) **@{-}, (9,6); (10,5) **@{-},
(5,0); (6,0) **@{-}, (6,1); (5,2) **@{-}, (4,1); (4,2) **@{-},
(8,1); (9,0) **@{-}, (10,0); (10,1) **@{-}, (8,2); (9,2) **@{-},
(2,4)*{\bullet}, (1,4)*{\bullet},(2,5)*{\bullet},(1,6)*{\bullet},(0,5)*{\bullet},(0,6)*{\bullet},
(5,4)*{\bullet}, (6,4)*{\bullet},(6,5)*{\bullet},(5,6)*{\bullet},(4,5)*{\bullet},(4,6)*{\bullet},
(9,4)*{\bullet}, (10,4)*{\bullet},(10,5)*{\bullet},(9,6)*{\bullet},(8,5)*{\bullet},(8,6)*{\bullet},
(5,0)*{\bullet}, (6,0)*{\bullet},(6,1)*{\bullet},(5,2)*{\bullet},(4,2)*{\bullet},(4,1)*{\bullet},
(9,0)*{\bullet}, (10,0)*{\bullet},(10,1)*{\bullet},(9,2)*{\bullet},(8,1)*{\bullet},(8,2)*{\bullet},
(-.3,5)*{2}, (-.3,6.3)*{1}, (1,3.7)*{1}, (2.3,3.7)*{1}, (2.3,5)*{2}, (1,6.3)*{2}, (5,3.7)*{1}, (6.3,3.7)*{2}, (6.3,5)*{2}, (5,6.3)*{2}, (3.7,6.3)*{1}, (3.7,5)*{1}, (9,3.7)*{2}, (10.3,3.7)*{2}, (10.3,5)*{1}, (9,6.3)*{2}, (7.7,6.3)*{1}, (7.7,5)*{1}, (3.7,1)*{2}, (3.7,2.3)*{1}, (5,2.3)*{2}, (6.3,1)*{1}, (6.3,-.3)*{2}, (5,-.3)*{1}, (7.7,1)*{1}, (7.7,2.3)*{1}, (9,2.3)*{2}, (10.3,1)*{2}, (10.3,-.3)*{1}, (9,-.3)*{2}
\end{xy}
}
}
},
(7,3.7)*{}
\end{xy}
\]

Notice that the top three matchings are obtained from one another by rotation, as are the two matchings in the second row. On the other hand, the corresponding standard tableaux are related by promotion: 
\[ \tableau{{1}&{3}&{4} \\ {2}&{5}&{6}} \xrightarrow{p} \tableau{{1}&{2}&{3} \\ {4}&{5}&{6}} \xrightarrow{p} \tableau{{1}&{2}&{5} \\ {3}&{4}&{6}} \xrightarrow{p} \tableau{{1}&{3}&{4} \\ {2}&{5}&{6}},\]
and 
\[\tableau{{1}&{3}&{5} \\ {2}&{4}&{6}} \xrightarrow{p} \tableau{{1}&{2}&{4} \\ {3}&{5}&{6}} \xrightarrow{p} \tableau{{1}&{3}&{5} \\ {2}&{4}&{6}}.\] 

In fact, by examining Yamanouchi words, the following theorem is easy to verify.

\begin{theorem}
Let $M$ denote a noncrossing matching on $2n$ vertices, and let $T$ be the corresponding standard Young tableau of shape $(n,n)$. Let $M'$ denote the noncrossing matching obtained by rotating $M$ clockwise by $\pi/n$. Then $p(T)$ is the tableau for $M'$.
\end{theorem}

This allows one in particular to deduce, or rather to see with one's own eyes, that promotion on a $2$ by $n$ rectangle has order dividing $2n$. It is natural to ask if such an elegant visualization of promotion is possible for other rectangles. We answer affirmatively for the three row case.

In \cite{K} Kuperberg introduced \emph{combinatorial rank $2$ spiders}. These are planar categories describing the invariant space $Inv(V_1 \otimes V_2 \otimes \dotsc \otimes V_n)$ of a tensor product of irreducible representations $V_i$ of a rank $2$ Lie algebra $\mathfrak g$. Spiders generalize the \emph{Temperley-Lieb category} that gives a similar basis for invariants of ${\mathfrak {sl}}_2$, see \cite{FKh}. Spiders are defined on a \emph{web space}: a vector space whose basis is a collection of planar graphs called \emph{webs}. These are the graphs we are interested in. In fact, the noncrossing matchings used above are exactly the $A_1$ webs, cf. \cite{K}. The question of describing spiders in arbitrary rank remains open. In this paper we restrict ourselvs to $A_1$ and $A_2$ spiders, and correspondingly to the two and three row cases.

The final part of our work deals with the \emph{cyclic sieving phenomenon} (CSP). Suppose we are given a finite set $X$, a finite cyclic group $C = \langle c \rangle$ acting on $X$, and a polynomial $X(q) \in \mathbb{Z}[q]$ with integer coefficients.  Following Reiner, Stanton, and White, \cite{RSW} we say that the triple $(X, C, X(q))$ exhibits the cyclic sieving phenomenon if for every integer $d \geq 0$, we have that $|X^{c^d}| = X(\zeta^d)$ where $\zeta \in \mathbb{C}$ is a root of unity of multiplicitive order $|C|$ and $X^{c^d}$ is the fixed point set of the action of the power $c^d$.  In particular, since the identity element fixes everything in any group action, we have that $|X| = X(1)$ whenever $(X, C, X(q))$ exhibits the CSP.  

If the triple $(X, C, X(q))$ exhibits the CSP and $\zeta$ is a primitive $|C|^{th}$ root of unity, we can determine the cardinalities of the fixed point sets $X^1 = X$, $X^c$, $X^{c^2}$, $\dots, X^{c^{|C|-1}}$ via the polynomial evaluations $X(1), X(\zeta), X(\zeta^2), \dots, X(\zeta^{|C|-1})$.  These fixed point set sizes determine the cycle structure of the canonical image of $c$ in the group of permutations of $X$, $S_X$.  Therefore, to find the cycle structure of the image of any bijection $c: X \rightarrow X$, it is enough to determine the order of the action of $c$ on $X$ and find a polynomial $X(q)$ such that $(X, \langle c \rangle, X(q))$ exhibits the CSP. 

In \cite{Rh} the third author proved an instance of the CSP related to the action of promotion on rectangular tableaux.  Recall that for any partition $\lambda \vdash n$, we have that the standard tableaux of shape $\lambda$ are enumerated by the Frame-Robinson-Thrall \emph{hook length formula}:
\[
|SYT(\lambda)| = \frac{n!}{\Pi_{(i,j) \in \lambda} h_{ij}},
\]
where the product is over the boxes $(i,j)$ in $\lambda$ and $h_{ij}$ is the \emph{hook length} at the box $(i,j)$, i.e., the number of boxes directly east or south of the box $(i,j)$ in $\lambda$, counting itself exactly once.  To obtain the polynomial used for cyclic sieving, we replace the hook length formula with a natural $q$-analogue. First, recall that for any $n \in \mathbb{N}$, $[n]_q := 1 + q + \cdots + q^{n-1}$ and $[n]_q! := [n]_q [n-1]_q \cdots [1]_q$.

\begin{theorem}[\cite{Rh}, Theorem 3.9] \label{thm:cs}
Let $\lambda \vdash n$ be a rectangular shape and let $X = SYT(\lambda)$.  Let $C := \mathbb{Z} / n \mathbb{Z}$ act on $X$ via promotion.  Then, the triple $(X, C, X(q))$ exhibits the cyclic sieving phenomenon, where
\begin{equation*}
X(q) = \frac{[n]_q!}{\Pi_{(i,j) \in \lambda} [h_{ij}]_q}
\end{equation*}
is the $q$-analogue of the hook length formula.
\end{theorem}

The proof in \cite{Rh} involves showing that the image of the long cycle $(n, n-1, \dots, 1) \in S_n$ in the Kazhdan-Lusztig cellular representation of shape $\lambda$ is, up to a predictable sign, the permutation matrix corresponding to the action of promotion on $SYT(\lambda)$, hence reducing the problem to a character evaluation.  This approach, while conceptually clean, has the drawback that it involves an object which is somewhat difficult to compute with and visualize---the KL cellular representation for rectangular shapes. Here we use webs as a basis for irreducible representations to give a simpler representation theoretic proof for the special cases of Theorem \ref{thm:cs} in which $\lambda$ has $2$ or $3$ rows.

The paper is structured as follows. In section \ref{sec:A2} we present $A_2$ webs along with some of their important known properties. We also state our first main result (Theorem \ref{thm:prom}), that rotation of $A_2$ webs is equivalent to promotion of rectangular tableaux with three rows. Because of its length, the proof of Theorem \ref{thm:prom} is relegated to Section \ref{sec:proof}. In Section \ref{sec:sieving} we give a self-contained proof of the cyclic sieving phenomenon for webs and derive some enumerative corollaries about rotational symmetry of webs. Section \ref{sec:concl} provides some ideas for future study.

\section{$A_2$-webs}\label{sec:A2}

Following Kuperberg \cite{K}, let us now define $A_2$ webs.

\begin{definition}
A planar directed graph $D$ with no multiple edges embedded in a disk is an \emph{$A_2$-web} if it satisfies the following conditions:
\begin{enumerate}
 \item $D$ is bipartite, with each edge of $D$ is oriented from one of the \emph{negative} vertices to one of the \emph{positive} vertices, and
 \item all the boundary vertices have degrees $1$ while all internal vertices have degree $3$.
\end{enumerate}
If, in addition, $D$ is \emph{non-elliptic}, i.e.,
\begin{enumerate}
 \item[(3)] all internal faces of $D$ have at least $6$ sides,
\end{enumerate}
then we say $D$ is an \emph{irreducible $A_2$ web}.
\end{definition}

When speaking of webs, we will omit the word irreducible when it is implied by the context. Note that webs are planar embeddings of graphs viewed up to a homeomorphism on the interior of the disk, with boundary vertices placed canonically.

Let $W^{(3)}$ denote the $\mathbb{C}$-vector space with basis the set of all irreducible $A_2$ webs.  Kuperberg \cite{K} introduced the following set of \emph{spider reduction rules} for $A_2$ webs.

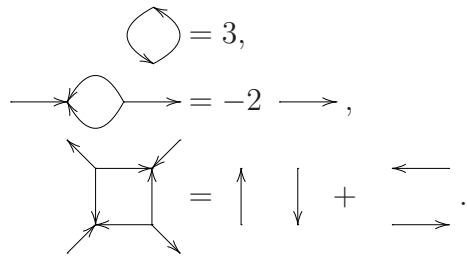
\begin{align*}
\begin{xy}
0;<.75cm,0cm>:
{\ar @{->} @/_10pt/ (0,-.5); (0,.5)}, {\ar @{<-} @/^10pt/ (0,-.5); (0,.5)}
\end{xy} &= 3, \\
\begin{xy}
0;<.75cm,0cm>:
{\ar @{->} (0,0); (1,0)}, {\ar @{<-} @/^10pt/ (1,0); (2,0)}, {\ar @{<-} @/_10pt/ (1,0); (2,0)}, {\ar @{->} (2,0); (3,0)}
\end{xy} &=-2 \;\; \begin{xy}
0;<.75cm,0cm>:
{\ar @{->} (0,0); (1,0)}
\end{xy} \;, \\
\begin{xy}
0;<.75cm,0cm>:
{\ar @{->} (0,-1); (.5,-.5)}, {\ar @{<-} (.5,-.5); (.5,.5)}, {\ar @{<-} (.5,-.5); (1.5,-.5)}, {\ar @{->} (.5,.5); (0,1)}, {\ar @{->} (1.5,-.5); (2,-1)}, 
{\ar @{->} (1.5,-.5); (1.5,.5)}, {\ar @{->} (.5,.5); (1.5,.5)}, {\ar @{<-} (1.5,.5); (2,1)} 
\end{xy} &= \; \,\, \begin{xy}
0;<.75cm,0cm>:
{\ar @{->} (0,-.5); (0,.5)}, {\ar @{->} (1,.5); (1,-.5)}
\end{xy} \,\, \;\; + \;\; \,\, \begin{xy}
0;<.75cm,0cm>:
{\ar @{->} (1,.5); (0,.5)}, {\ar @{->} (0,-.5); (1,-.5)}
\end{xy}\,\,.
\end{align*}

These local graph transformations, when iterated, allow for the expression of an arbitrary $A_2$ web as a linear combination of irreducible $A_2$ webs.  Moreover, it can be shown \cite{K2} that any application of these rules to a fixed $A_2$ web yields the same linear combination of irreducible $A_2$ webs. In other words the spider reduction rules are {\it {confluent}}.

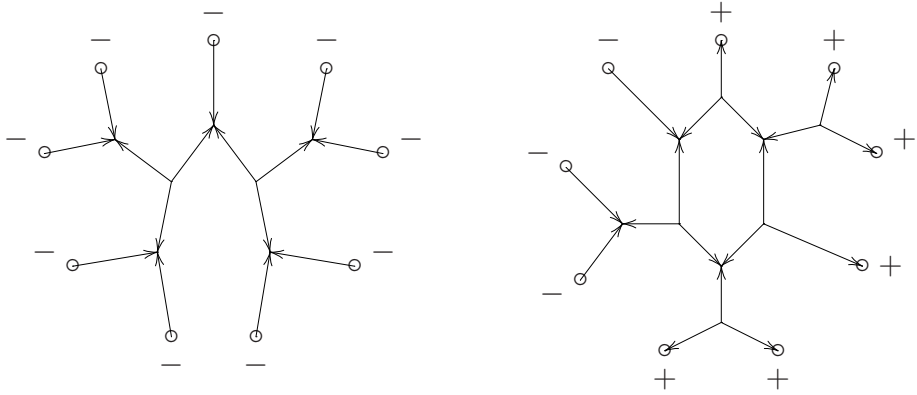
\begin{figure}[h!]
\begin{center}
\[
\begin{xy}
0;<.75cm,0cm>:
(3,6)="l1", (3,6.5)*{-}, (5,5.5)="l2", (5,6)*{-}, (6,4)="l3", (6.5,4.25)*{-}, (5.5,2)="l4", (6,2.25)*{-}, (3.75,.75)="l5", (3.75,.25)*{-}, (2.25,.75)="l6", (2.25,.25)*{-}, (.5,2)="l7", (0,2.25)*{-}, (0,4)="l8", (-.5,4.25)*{-}, (1,5.5)="l9", (1,6)*{-}, "l1"*{\circ}, "l2"*{\circ}, "l3"*{\circ}, "l4"*{\circ}, "l5"*{\circ}, "l6"*{\circ}, "l7"*{\circ}, "l8"*{\circ}, "l9"*{\circ},
(3,4.5)="a1", (3.75,3.5)="a2", (4.75,4.25)="a3", (2.25,3.5)="a4", (1.25,4.25)="a5", (2,2.25)="a6", (4,2.25)="a7",
{\ar @{->} "l1"; "a1"}, {\ar @{->} "l2"; "a3"},{\ar @{->} "l3"; "a3"}, {\ar @{->} "l4"; "a7"}, {\ar @{->} "l5"; "a7"}, {\ar @{->} "l6"; "a6"}, {\ar @{->} "l7"; "a6"}, {\ar @{->} "l8"; "a5"}, {\ar @{->} "l9"; "a5"}, {\ar @{->} "a2"; "a1"}, {\ar @{->} "a2"; "a3"}, {\ar @{->} "a2"; "a7"}, {\ar @{->} "a4"; "a1"}, {\ar @{->} "a4"; "a5"}, {\ar @{->} "a4"; "a6"},
(12,6)="r1", (12,6.5)*{+}, (14,5.5)="r2", (14,6)*{+}, (14.75,4)="r3", (15.25,4.25)*{+}, (14.5,2)="r4", (15,2)*{+}, (13,.5)="r5", (13,0)*{+}, (11,.5)="r6", (11,0)*{+}, (9.5,1.75)="r7", (9,1.5)*{-},(9.25,3.75)="r8", (8.75,4)*{-}, (10,5.5)="r9", (10,6)*{-}, "r1"*{\circ}, "r2"*{\circ}, "r3"*{\circ}, "r4"*{\circ}, "r5"*{\circ}, "r6"*{\circ}, "r7"*{\circ}, "r8"*{\circ}, "r9"*{\circ},
(12,5)="b1", (12.75,4.25)="b2", (12.75,2.75)="b3", (12,2)="b4", (11.25,2.75)="b5", (11.25,4.25)="b6", (10.25,2.75)="b7", (12,1)="b8", (13.75,4.5)="b9",
{\ar @{<-} "r1"; "b1"}, {\ar @{<-} "r2"; "b9"},{\ar @{<-} "r3"; "b9"},{\ar @{<-} "r4"; "b3"},{\ar @{<-} "r5"; "b8"},{\ar @{<-} "r6"; "b8"}, {\ar @{->} "r7";"b7"}, {\ar @{->} "r8";"b7"}, {\ar @{->} "r9";"b6"}, {\ar @{->} "b5";"b7"}, {\ar @{->} "b5";"b6"}, {\ar @{->} "b5";"b4"}, {\ar @{<-} "b2";"b9"}, {\ar @{<-} "b2";"b1"}, {\ar @{<-} "b2";"b3"}, {\ar @{->} "b1";"b6"}, {\ar @{<-} "b4";"b8"}, {\ar @{<-} "b4";"b3"}
\end{xy}
\]
\end{center}
\caption{Two irreducible $A_2$ webs.}\label{fig:pr1}
\end{figure}

Figure \ref{fig:pr1} shows two examples of irreducible webs; the signs of boundary vertices are marked. Let $\gamma$ denote a cyclically ordered arrangement of signs (i.e., $+$ or $-$). We write $|\gamma|=n$ if the total number of boundary vertices is $n$. Let $\mathfrak M_{\gamma}$ denote the set of all irreducible webs with boundary $\gamma$. The following is a specialization of a theorem of Kuperberg.

\begin{theorem}[\cite{K}, Theorem 6.1] \label{thm:eq}
 Let $\gamma$ be a fixed boundary with $k$ ``$+$''s and $3n-2k$ ``$-$"s. The number of semistandard Young tableaux of shape $(3,\ldots, 3)$ and type $\{1^2,\ldots,k^2,k+1,\ldots,3n-k\}$ is equal to the cardinality of $\mathfrak M_{\gamma}$. 
\end{theorem}

In particular, if $\gamma$ has $3n$ ``$-$"s, i.e., $k=0$, then the set $\mathfrak M_n := \mathfrak M_{\gamma}$ and $SYT((n,n,n))$ are equinumerous. Kuperberg and Khovanov \cite{KKh} give an explicit bijection between these two sets. We now describe this bijection.

Place the boundary vertices of a web $D \in \mathfrak M_n$ on a line so that the web hangs below the line. We need to make a choice here where to cut the circular boundary. Next, consider the set of faces $F(D)$ created by the web and the line. Distinguish the infinite {\it {outer face}} $f_0$. For each $f \in F(D)$ we let the {\it {depth}}  of $f$, $d(f)$, be the minimal number of edges in $D$ one needs to cross to reach $f_0$ starting in $f$. In particular, $d(f_0)=0$. For an edge $e$ of $D$ let $f_e^l$ and $f_e^r$ denote the faces to the left and to the right of $e$ looking in the direction of $e$'s orientation. Label each edge $e$ of $D$ with the label $l(e) = d(f_e^l)-d(f_e^r)$.

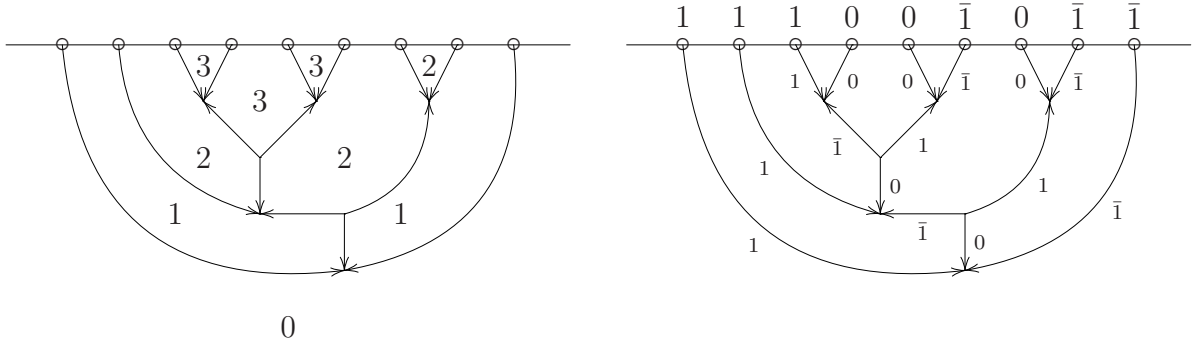
\begin{figure}[h!]
\begin{center}
\[
\begin{xy}
0;<.75cm,0cm>:
(0,5); (10,5) **@{-},
(1,5)="l1", (2,5)="l2", (3,5)="l3", (4,5)="l4",  (5,5)="l5", (6,5)="l6", (7,5)="l7", (8,5)="l8", (9,5)="l9", "l1"*{\circ}, "l2"*{\circ}, "l3"*{\circ}, "l4"*{\circ}, "l5"*{\circ}, "l6"*{\circ}, "l7"*{\circ}, "l8"*{\circ}, "l9"*{\circ},
(3.5,4)="a1", (5.5,4)="a2", (7.5,4)="a3", (4.5,3)="a4", (4.5,2)="a5", (6,2)="a6", (6,1)="a7",
{\ar @{->} @/_35pt/ "l1"; "a7"}, {\ar @{->} @/_15pt/ "l2"; "a5"}, {\ar @{->} "l3"; "a1"}, {\ar @{->} "l4"; "a1"}, {\ar @{->} "l5"; "a2"}, {\ar @{->} "l6"; "a2"}, {\ar @{->} "l7"; "a3"}, {\ar @{->} "l8"; "a3"}, {\ar @{->} @/^25pt/ "l9"; "a7"}, {\ar @{->} "a4"; "a1"}, {\ar @{->} "a4"; "a2"}, {\ar @{->} "a4"; "a5"}, {\ar @{->} "a6"; "a5"}, {\ar @{->} @/_10pt/ "a6"; "a3"}, {\ar @{->} "a6"; "a7"},
(5,0)*{0}, (3,2)*{1}, (7,2)*{1}, (6,3)*{2}, (3.5,3)*{2}, (4.5,4)*{3}, (3.5,4.625)*{3}, (5.5,4.625)*{3}, (7.5,4.625)*{2},
(11,5); (21,5) **@{-},
(12,5)="r1", (12,5.5)*{1}, (13,5)="r2", (13,5.5)*{1}, (14,5)="r3", (14,5.5)*{1}, (15,5)="r4", (15,5.5)*{0}, (16,5)="r5", (16,5.5)*{0}, (17,5)="r6", (17,5.5)*{\bar{1}}, (18,5)="r7", (18,5.5)*{0}, (19,5)="r8", (19,5.5)*{\bar{1}}, (20,5)="r9", (20,5.5)*{\bar{1}}, "r1"*{\circ}, "r2"*{\circ}, "r3"*{\circ}, "r4"*{\circ}, "r5"*{\circ}, "r6"*{\circ}, "r7"*{\circ}, "r8"*{\circ}, "r9"*{\circ},
(14.5,4)="b1", (16.5,4)="b2", (18.5,4)="b3", (15.5,3)="b4", (15.5,2)="b5", (17,2)="b6", (17,1)="b7",
{\ar @{->}_1 @/_35pt/ "r1"; "b7"}, {\ar @{->}_1 @/_15pt/ "r2"; "b5"}, {\ar @{->}_1 "r3"; "b1"}, {\ar @{->}^0 "r4"; "b1"}, {\ar @{->}_0 "r5"; "b2"}, {\ar @{->}^{\bar{1}} "r6"; "b2"}, {\ar @{->}_0 "r7"; "b3"}, {\ar @{->}^{\bar{1}} "r8"; "b3"}, {\ar @{->}^{\bar{1}} @/^25pt/ "r9"; "b7"}, {\ar @{->}^{\bar{1}} "b4"; "b1"}, {\ar @{->}_1 "b4"; "b2"}, {\ar @{->}^0 "b4"; "b5"}, {\ar @{->}^{\bar{1}} "b6"; "b5"}, {\ar @{->}_1 @/_10pt/ "b6"; "b3"}, {\ar @{->}^0 "b6"; "b7"},
\end{xy}
\]
\end{center}
\caption{Depths of faces and edge labelings.}\label{fig:pr2}
\end{figure}

Using the web on the left of Figure \ref{fig:pr1} as an example, we see it stretched out and labeled in Figure \ref{fig:pr2}. Note that the depth of two adjacent faces differs by at most $1$, which implies that each edge label is either $-1$, $0$ or $1$. In particular, one can read off the sequence of labels assigned to boundary edges. The web on Figure \ref{fig:pr2} thus produces the sequence $(1,1,1,0,0,-1,0,-1-1)$, which we can also write as a word $w = w(D) = 11100\bar{1}0\bar{1}\bar{1}$ (with $\bar{1}$ for $-1$). Any such resulting word $w(D)$ is \emph{dominant} in the language of \cite{KKh}, see \cite[Proposition 1]{KKh} and the preceding discussion. In our terminology, this means it is a balanced Yamanouchi word on the multiset $\{1^n,0^n,\bar{1}^n\}$. As mentioned earlier, such words are in bijection with standard Young tableaux of shape $(n,n,n)$, (here $1$ corresponds to row $1$, $0$ corresponds to row $2$, and $\bar{1}$ to row $3$).

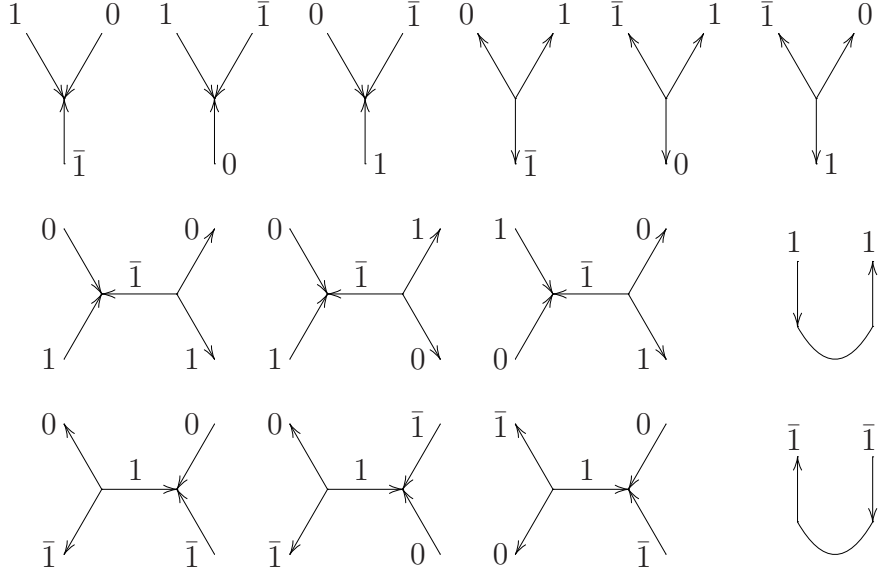
\begin{figure}[h!]
\begin{center}
\begin{xy}
0;<1cm,0cm>:<.5cm,\halfrootthree cm>::
(.5,0)*{};
{\ar @{->} (0,5);(1,4)}, {\ar @{->} (1,5);(1,4)}, {\ar @{->} (1.5,3);(1,4)},
{\ar @{->} (2,5);(3,4)}, {\ar @{->} (3,5);(3,4)}, {\ar @{->} (3.5,3);(3,4)},
{\ar @{->} (4,5);(5,4)}, {\ar @{->} (5,5);(5,4)}, {\ar @{->} (5.5,3);(5,4)},
{\ar @{<-} (6,5);(7,4)}, {\ar @{<-} (7,5);(7,4)}, {\ar @{<-} (7.5,3);(7,4)},
{\ar @{<-} (8,5);(9,4)}, {\ar @{<-} (9,5);(9,4)}, {\ar @{<-} (9.5,3);(9,4)}, 
{\ar @{<-} (10,5);(11,4)}, {\ar @{<-} (11,5);(11,4)}, {\ar @{<-} (11.5,3);(11,4)}, 
{\ar @{->} (14,-2.5);(13.5,-1.5)}, {\ar @{<-} (15,-2.5);(14.5,-1.5)}, (14,-2.5);(15,-2.5) **\crv{(15,-3.5)},
{\ar @{->} (12,1.5);(12.5,.5)}, {\ar @{<-} (13,1.5);(13.5,.5)}, (12.5,.5); (13.5,.5) **\crv{(13.5,-.5)},
{\ar @{->} (2,2);(3,1)}, {\ar @{->} (3,0);(3,1)}, {\ar @{<-} (3,1);(4,1)}, {\ar @{->} (4,1);(4,2)}, {\ar @{->} (4,1);(5,0)},
{\ar @{->} (5,2);(6,1)}, {\ar @{->} (6,0);(6,1)}, {\ar @{<-} (6,1);(7,1)}, {\ar @{->} (7,1);(7,2)}, {\ar @{->} (7,1);(8,0)},
{\ar @{->} (8,2);(9,1)}, {\ar @{->} (9,0);(9,1)}, {\ar @{<-} (9,1);(10,1)}, {\ar @{->} (10,1);(10,2)}, {\ar @{->} (10,1);(11,0)},
{\ar @{<-} (3.5,-1);(4.5,-2)}, {\ar @{<-} (4.5,-3);(4.5,-2)}, {\ar @{->} (4.5,-2);(5.5,-2)}, {\ar @{<-} (5.5,-2);(5.5,-1)}, {\ar @{<-} (5.5,-2);(6.5,-3)},
{\ar @{<-} (6.5,-1);(7.5,-2)}, {\ar @{<-} (7.5,-3);(7.5,-2)}, {\ar @{->} (7.5,-2);(8.5,-2)}, {\ar @{<-} (8.5,-2);(8.5,-1)}, {\ar @{<-} (8.5,-2);(9.5,-3)},
{\ar @{<-} (9.5,-1);(10.5,-2)}, {\ar @{<-} (10.5,-3);(10.5,-2)}, {\ar @{->} (10.5,-2);(11.5,-2)}, {\ar @{<-} (11.5,-2);(11.5,-1)}, {\ar @{<-} (11.5,-2);(12.5,-3)},
(-.3,5.3)*{1}, (1,5.3)*{0}, (1.7,3)*{\bar{1}},
(1.7,5.3)*{1}, (3,5.3)*{\bar{1}}, (3.7,3)*{0},
(3.7,5.3)*{0}, (5,5.3)*{\bar{1}}, (5.7,3)*{1},
(5.7,5.3)*{0}, (7,5.3)*{1}, (7.7,3)*{\bar{1}},
(7.7,5.3)*{\bar{1}}, (9,5.3)*{1}, (9.7,3)*{0},
(9.7,5.3)*{\bar{1}}, (11,5.3)*{0}, (11.7,3)*{1},
(11.8,1.8)*{1}, (12.8,1.8)*{1},
(13.3,-1.2)*{\bar{1}}, (14.3,-1.2)*{\bar{1}},
(1.8,2)*{0}, (3.7,2)*{0}, (3.3,1.3)*{\bar{1}}, (2.8,0)*{1}, (4.7,0)*{1},
(4.8,2)*{0}, (6.7,2)*{1}, (6.3,1.3)*{\bar{1}}, (5.8,0)*{1}, (7.7,0)*{0},
(7.8,2)*{1}, (9.7,2)*{0}, (9.3,1.3)*{\bar{1}}, (8.8,0)*{0}, (10.7,0)*{1},
(3.3,-1)*{0}, (5.2,-1)*{0}, (4.8,-1.7)*{1}, (4.3,-3)*{\bar{1}}, (6.2,-3)*{\bar{1}},
(6.3,-1)*{0}, (8.2,-1)*{\bar{1}}, (7.8,-1.7)*{1}, (7.3,-3)*{\bar{1}}, (9.2,-3)*{0},
(9.3,-1)*{\bar{1}}, (11.2,-1)*{0}, (10.8,-1.7)*{1}, (10.3,-3)*{0}, (12.2,-3)*{\bar{1}}
\end{xy}
\end{center}
\caption{Growth rules for labeled $A_2$ webs.}\label{fig:pr3}
\end{figure}

In order to define the inverse map, that is, how to assign a unique web $D(w)$ to every dominant word $w$, we need the {\it {growth rules}} given in Figure \ref{fig:pr3}. These pictures describe local moves for joining together dangling ``strands" according to their orientation and labeling, and can be used to generate any irreducible web. Given a sign sequence $\gamma$ and a word $w$ with three distinct letters $1 < 0 < \bar{1}$, we first draw vertices on a line, labeled from left to right by $w$. Then we draw a directed edge downward from each vertex. To form the web, we choose a pair of neighboring strands (i.e., with no strands dangling between) and apply the local rules in Figure \ref{fig:pr3} to join the strands together. We continue until there are no neighboring strands to which we can apply the  growth rules.

\begin{remark}
 The growth rules here are a slight modification of growth rules in \cite{KKh}, but are nonetheless equivalent. We have defined our rules so that the induced edge labelings are consistent with the depths of the faces of $D$. See Lemma \ref{lem:1}. To obtain the Khovanov-Kuperberg rules from ours, ignore all horizontal labels and negate the labels for upward pointing arrows. It is straightforward (if tedious) to verify that our modified rules give rise to the same claims asserted in Lemmas 1--3, and Proposition 1 of \cite{KKh}, summarized in Theorem \ref{thm:kkh} below.
\end{remark}

The following is the compilation of several statements proved by Khovanov and Kuperberg.

\begin{theorem}[\cite{KKh}, Lemmas 1--3, Proposition 1] \label{thm:kkh}
The web produced by the growth algorithm does not depend on the choices made in applying the growth rules. Furthermore, if one starts with a dominant word $w$ and a sign sequence of all ``$+$"s or all ``$-$"s, the growth algorithm does not terminate until there are no dangling strands, and when it terminates the resulting web is non-elliptic. In fact, the maps $w$ and $D$ are inverses in this case and provide a bijection between irreducible webs and dominant (i.e., balanced Yamanouchi) words.
\end{theorem}

Finally, we are ready to state and prove the result relating webs with promotion. Let $p(D)$ be the web obtained by rotating a web $D$ by $\frac{2\pi}{3n}$, so that if we cut the boundary in the same place, the first vertex on the boundary of $D$ becomes the last vertex on the boundary of $p(D)$.

\begin{theorem} \label{thm:prom}
For $D\in \mathfrak{M}_n$, we have \[T(w(p(D))) = p(T(w(D))).\] That is, the tableau associated with the rotation of $D$ is given by promotion of the tableau associated with $D$ itself.
\end{theorem}

\begin{figure}[h!]
\begin{center}
\[
\begin{xy}
0;<.75cm,0cm>:
{\ar @{->}^p (9.75,3); (11.25,3)},
(0,5); (10,5) **@{-},
(1,5)="l1", (1,5.5)*{1}, (2,5)="l2", (2,5.5)*{1}, (3,5)="l3", (3,5.5)*{1}, (4,5)="l4", (4,5.5)*{0},  (5,5)="l5", (5,5.5)*{0}, (6,5)="l6", (6,5.5)*{\bar{1}}, (7,5)="l7", (7,5.5)*{0}, (8,5)="l8", (8,5.5)*{\bar{1}}, (9,5)="l9", (9,5.5)*{\bar{1}}, "l1"*{\circ}, "l2"*{\circ}, "l3"*{\circ}, "l4"*{\circ}, "l5"*{\circ}, "l6"*{\circ}, "l7"*{\circ}, "l8"*{\circ}, "l9"*{\circ},
(3.5,4)="a1", (5.5,4)="a2", (7.5,4)="a3", (4.5,3)="a4", (4.5,2)="a5", (6,2)="a6", (6,1)="a7",
{\ar @{->} @/_35pt/ "l1"; "a7"}, {\ar @{->} @/_15pt/ "l2"; "a5"}, {\ar @{->} "l3"; "a1"}, {\ar @{->} "l4"; "a1"}, {\ar @{->} "l5"; "a2"}, {\ar @{->} "l6"; "a2"}, {\ar @{->} "l7"; "a3"}, {\ar @{->} "l8"; "a3"}, {\ar @{->} @/^25pt/ "l9"; "a7"}, {\ar @{->} "a4"; "a1"}, {\ar @{->} "a4"; "a2"}, {\ar @{->} "a4"; "a5"}, {\ar @{->} "a6"; "a5"}, {\ar @{->} @/_10pt/ "a6"; "a3"}, {\ar @{->} "a6"; "a7"},
(11,5); (21,5) **@{-},
(12,5)="r1", (12,5.5)*{1}, (13,5)="r2", (13,5.5)*{1}, (14,5)="r3", (14,5.5)*{0}, (15,5)="r4", (15,5.5)*{0}, (16,5)="r5", (16,5.5)*{\bar{1}}, (17,5)="r6", (17,5.5)*{1}, (18,5)="r7", (18,5.5)*{\bar{1}}, (19,5)="r8", (19,5.5)*{0}, (20,5)="r9", (20,5.5)*{\bar{1}}, "r1"*{\circ}, "r2"*{\circ}, "r3"*{\circ}, "r4"*{\circ}, "r5"*{\circ}, "r6"*{\circ}, "r7"*{\circ}, "r8"*{\circ}, "r9"*{\circ},
(13.5,4)="b1", (15.5,4)="b2", (17.5,4)="b3", (19.5,4)="b4", (14.5,3)="b5", (14.5,2)="b6", (17.5,2.5)="b7",
{\ar @{->} @/_20pt/ "r1"; "b6"}, {\ar @{->} "r2"; "b1"}, {\ar @{->} "r3"; "b1"}, {\ar @{->} "r4"; "b2"}, {\ar @{->} "r5"; "b2"}, {\ar @{->} "r6"; "b3"}, {\ar @{->} "r7"; "b3"}, {\ar @{->} "r8"; "b4"}, {\ar @{->} "r9"; "b4"}, {\ar @{->} "b5"; "b1"}, {\ar @{->} "b5"; "b2"}, {\ar @{->} "b5"; "b6"}, {\ar @{->} @/^5pt/ "b7"; "b6"}, {\ar @{->} "b7"; "b3"}, {\ar @{->} @/_10pt/ "b7"; "b4"}
\end{xy}
\]
\end{center}
\caption{Rotation of an $A_2$ web.}\label{fig:pr5}
\end{figure}
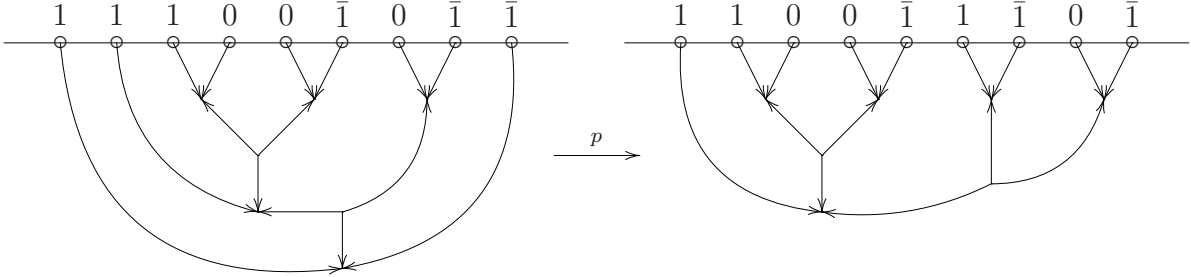

\begin{example}
Figure \ref{fig:pr5} shows an example of cyclic rotation of a web. The corresponding map on standard Young tableaux is:
$$
\tableau{{1}&{2}&{3} \\ {4}&{5}&{7} \\ {6}&{8}&{9}} \quad
\xrightarrow{p} \quad
\tableau{{1}&{2}&{6} \\ {3}&{4}&{8} \\ {5}&{7}&{9}} \,.
$$
\end{example}

\section{Proof of Theorem \ref{thm:prom}}\label{sec:proof}

Throughout this section we assume $D$ is irreducible.

The main idea for the proof of Theorem \ref{thm:prom} is as follows. Given a web $D$ with word $w$, we describe a way to cut it into three regions: $\Dl$, $\Dm$, and $\Dr$; see Figure \ref{fig:pr8}. Upon moving the leftmost vertex of $D$ to the right, the depths of faces in $\Dl$ decrease by $1$, the depths of faces in $\Dr$ increase by $1$, and the depths of faces in $\Dm$ remain unchanged. In terms of the word obtained after rotation, $w' = w(p(D))$, this means that we have \[w' = w_2 \cdots w_{a-1} 1 w_{a+1} \cdots w_{b-1} 0 w_{b+1} \cdots w_{3n} \bar{1},\] where $a$ is the position of the boundary vertex on the border between $\Dl$ and $\Dm$ and the boundary vertex in position $b$ lies between $\Dm$ and $\Dr$. We want to show that $w' = p(w)$, which is to say that $a$ is the first position where letters $1$ and $0$ balance in $w$, and $b$ is the first position after $a$ in $w$ where letters $0$ and $\bar{1}$ balance. 

We start with proving some lemmas. For a web in $\mathfrak M_n$ let us say that two consecutive boundary vertices are {\it {neighbors}} if they are connected to a common internal vertex. For example, the web in Figure \ref{fig:pr2} has four pairs of neighbors. 

\begin{lemma} \label{lem:3}
 Every web in $\mathfrak M_n$  has at least three pairs of neighbors.
\end{lemma}

\begin{proof}
 Let $D \in \mathfrak M_n$ be a web and let $w = w(D)$ be the associated dominant word. We are looking for occurances of $\cdots 10 \cdots$, $\cdots 1\bar{1} \cdots$ and $\cdots 0\bar{1} \cdots$ in $w$. It is clear that there are at least two such pairs: at the first appearance of a $0$, which is preceded immediately by a $1$, and at the final appearance of a $0$, immediately followed by a $\bar{1}$. Thus, the lemma follows if there are two or more connected components.

If \[w = \underbrace{1\cdots 1}_n\underbrace{0\cdots 0}_n\underbrace{\bar{1}\cdots \bar{1}}_n,\] then the first $1$ and the final $\bar{1}$ also form a pair of neighbors. (An easy induction argument using the growth rules shows that the word \[ \underbrace{1\cdots 1}_n\underbrace{0\cdots 0}_n\] has $n$ upward-pointing strands dangling, each labeled with $\bar{1}$, and the leftmost of these is adjacent to the edge from $w_1$. These pair off with the remaining letters $\bar{1}$ in $w$ without intersecting, forcing the edges from $w_1$ and $w_n$ to be adjacent.) 

If $w$ is not of this form, then either 
\begin{enumerate}
 \item the last $1$ occurs after the first $0$, in which case there is an extra occurrence of $\cdots 10 \cdots$ or a $\cdots 1\bar{1} \cdots$, or
 \item the first $\bar{1}$ occurs before the last $0$, in which case there is an occurrence of a $\cdots 1\bar{1} \cdots$ or an extra occurrence of $\cdots 0\bar{1} \cdots$.
\end{enumerate}
\end{proof}

Let us now take an edge $e$ adjacent to a vertex $v$ of an irreducible web $D$, and define a {\it {left cut}} \[C^l_{e,v}: \xrightarrow{e} v \xrightarrow{e_1^l} v_1^l \xrightarrow{e_2^l} \cdots \xrightarrow{e_i^l} v_i^l \xrightarrow{e_{i+1}^l} \cdots \] and a {\it right cut} \[C^r_{e,v}: \xrightarrow{e} v \xrightarrow{e_1^r} v_1^r \xrightarrow{e_2^r} \cdots \xrightarrow{e_j^r} v_j^r \xrightarrow{e_{j+1}^r} \cdots \] starting at $v$ as follows. We move along $e$ towards $v$ (the orientation of $e$ does not matter here) and turn left at $v$ onto edge $e_1^l$. Traversing $e_1^l$ we reach the next vertex $v_1^l$ where we turn right onto edge $e_2^l$. This takes us to vertex $v_2^l$, and so on. We keep alternating left and right turns until we reach a boundary vertex, at which point the process stops. The left cut $C^l_{e,v}$ is the resulting sequence of edges and vertices. Similarly we define the right cut with edges $e_j^r$ and vertices $v_j^r$, with the only difference being that the first turn at $v$ is to the right.

\begin{lemma} \label{lem:4}
 For any $e$ and $v$ the left and right cuts do not intersect each other and do not self-intersect. In other words, all vertices $v, v_i^l, v_j^r$ are distinct.
\end{lemma}

\begin{proof}
Recall that all internal faces of an irreducible web must have at least $6$ sides. We will show that if the left and right cuts intersect (or self-intersect) then the web must have a $4$-cycle, a contradiction.

Let $D$ be a web, and consider the left and right cut for a given pair $(e,v)$. Assume that the left cut intersects the right cut, and take the first point of intersection, $w$. There are several possible scenarios to consider, based on the sign of $w$ and $v$ and on whether the third edge at $w$ points inward or outward with respect to the enclosed region. One of the cases is shown in Figure \ref{fig:pr7}. In fact, this scenario is in some sense the ``worst" one. The key observation is that the part of the original web contained inside the cycle formed by left and right cuts is also a web, say $D'$, with all the boundary edges of the same orientation. 

\begin{figure}[h!]
\begin{center}
\[
\begin{xy}
0;<.75cm,0cm>:
(3,1)="v", (3.25,6.25)="w", {\ar @{-}^e (3,-.5); "v"}, (3.25,7.25); "w" **@{-},
(2.5,1.5) ="l1", (2.25,2.25)="l2", (1.5,2.75)="l3", (1.75,3.5)="l4", (1.5,4.25)="l5", (2,5)="l6", (2.35,5.85)="l7", "v"; "l1" **@{-}, "l1"; "l2" **@{-}, (1.75,1) **@{-}, "l2"; "l3" **@{-}, (3,2.5) **@{-}, "l3"; "l4" **@{-}, (.75,2.5) **@{-}, "l4"; "l5" **@{-}, (2.5,3.45) **@{-}, "l5"; "l6" **@{-}, (.75, 4.5) **@{-}, "l6"; "l7" **@{-}, (2.8,4.75) **@{-}, "l7"; "w" **@{-}, (1.85,6.35) **@{-},
(4,1.1) ="r1", (4.5,1.75)="r2", (5.25,2.25)="r3",(4.75,3)="r4",(5.25,4)="r5",(4.55,4.8)="r6",(4.25,6)="r7",
"v"; "r1" **@{-}, "r1"; "r2" **@{-}, (4.25,.5) **@{-}, "r2"; "r3" **@{-}, (3.9,2.25) **@{-}, "r3"; "r4" **@{-}, (6,2.2) **@{-}, "r4"; "r5" **@{-}, (4,3) **@{-}, "r5"; "r6" **@{-}, (6,4) **@{-}, "r6"; "r7" **@{-}, (4, 4.3) **@{-}, "r7"; "w" **@{-}, (4.75,6.5) **@{-},
"v"*{\bullet}, "w"*{\bullet}, "l1"*{\circ}, "l2"*{\bullet}, "l3"*{\circ}, "l4"*{\bullet}, "l5"*{\circ}, "l6"*{\bullet}, "l7"*{\circ}, "r1"*{\circ}, "r2"*{\bullet}, "r3"*{\circ}, "r4"*{\bullet}, "r5"*{\circ}, "r6"*{\bullet}, "r7"*{\circ}, (3.25,1.25)*{v}, (3.25,5.85)*{w}, (3.25,3.5)*{D'}
\end{xy}
\]
\end{center}
\caption{Intersection of left and right cuts.}\label{fig:pr7}
\end{figure}
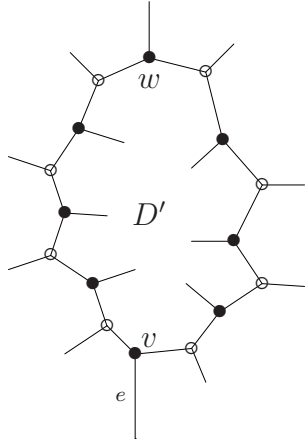

If the web $D'$ is empty then the left and right cut form a $4$-cycle and we are done. 

Assume that $D'$ is nonempty. Then by the final claim of Theorem \ref{thm:kkh} there are at least $3$ boundary vertices in $D'$. By Lemma \ref{lem:3}, $D'$ must have at least $3$ pairs of neighboring vertices. One of these pairs might be placed next to  vertex $v$ and another pair may sit next to vertex $w$, but the third one must occur somewhere in between, and as seen in Figure \ref{fig:pr7}, it unavoidably creates a $4$-cycle.

Other cases to consider are very similar, as well as the non-self-intersection claim.
\end{proof}

\begin{lemma} \label{lem:1}
Given a dominant word $w$, the edge labels given in the growth rules are consistent with the depths of the faces of the resulting web $D(w)$.
\end{lemma}

\begin{proof}
The proof is by induction on the maximal depth of a face. Clearly if the maximal depth is 1, then $w$ consists of copies of $10\bar{1}$, and the labeling is consistent.

Let $D$ be a web of depth at least 2 corresponding to a dominant word. Since it is dominant, Theorem \ref{thm:kkh} tells us there are no dangling strands. Further notice that the growth algorithm must finish with two $1$s of opposite parity connecting up or two $\bar{1}$s connecting up. In either of these situations the labeling is appropriate since we are creating a face of depth $1$. Moreover, since the confluence property allows us to perform the growth operations in any order, we see that the boundary between the inner faces and the outer face $f_0(D)$ consists of edges alternately labeled with $1$ and $\bar{1}$.

Remove from $D$ all of these edges and the $0$-labeled edges attached to them, which we call the \emph{outer strip} of $D$. Since these $0$-labeled edges separate faces of depth $1$, these labelings are also consistent with measurement of depth. In Figure \ref{fig:strip}, the outer strip is indicated with dashed lines.

\begin{figure}[h!]
\[
\begin{xy}
0;<.75cm,0cm>:
(11,5); (21,5) **@{-},
(12,5)="r1", (12,5.5)*{1}, (13,5)="r2", (13,5.5)*{1}, (14,5)="r3", (14,5.5)*{0}, (15,5)="r4", (15,5.5)*{0}, (16,5)="r5", (16,5.5)*{\bar{1}}, (17,5)="r6", (17,5.5)*{1}, (18,5)="r7", (18,5.5)*{\bar{1}}, (19,5)="r8", (19,5.5)*{0}, (20,5)="r9", (20,5.5)*{\bar{1}}, "r1"*{\circ}, "r2"*{\circ}, "r3"*{\circ}, "r4"*{\circ}, "r5"*{\circ}, "r6"*{\circ}, "r7"*{\circ}, "r8"*{\circ}, "r9"*{\circ},
(13.5,4)="b1", (15.5,4)="b2", (17.5,4)="b3", (19.5,4)="b4", (14.5,3)="b5", (14.5,2)="b6", (17.5,2.5)="b7",
{\ar @{-->}_1 @/_20pt/ "r1"; "b6"}, {\ar @{->}_1 "r2"; "b1"}, {\ar @{->}^0 "r3"; "b1"}, {\ar @{->}_0 "r4"; "b2"}, {\ar @{->}^{\bar{1}} "r5"; "b2"}, {\ar @{->}_1 "r6"; "b3"}, {\ar @{->}^{\bar{1}} "r7"; "b3"}, {\ar @{-->}_0 "r8"; "b4"}, {\ar @{-->}^{\bar{1}} "r9"; "b4"}, {\ar @{->}^{\bar{1}} "b5"; "b1"}, {\ar @{->}_1 "b5"; "b2"}, {\ar @{-->}^0 "b5"; "b6"}, {\ar @{-->}^{\bar{1}} @/^5pt/ "b7"; "b6"}, {\ar @{-->}_0 "b7"; "b3"}, {\ar @{-->}_1 @/_10pt/ "b7"; "b4"}
\end{xy}
\]
\caption{The outer strip of a web.}\label{fig:strip}
\end{figure}
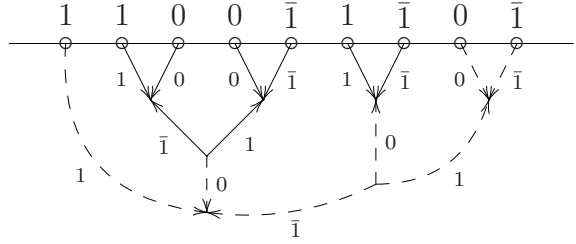

What remains, $D'$, is not necessarily a web, but it is a planar edge-labeled graph with the faces of $D$ of depth at least $2$. The graph $D'$ will have one or more connected components which we must examine individually.

Let $D''$ be one of these connected components. If it forms a web by itself, we are done by the induction hypothesis. Suppose $D''$ is a connected component which is not a proper web by itself, i.e., at least one of the $0$-labeled edges on the outer strip we removed from $D$ was connected to $D''$. We claim that its boundary with the outer face of $D''$ consists of edges labeled alternately $1$ and $\bar{1}$, at which point the proof again follows by induction, the base cases being of the form $w(D'') = 10\cdots0\bar{1}$, where all closed regions have depth 1 and the labeling is easily verified.

Let us consider approaching $D''$ along one of the $0$-labeled edges on the outer strip of $D$. By looking at the local growth rules, we see there are six possibilities for the neighborhood of the point where the $0$-labeled edge meets $D''$. In each case, the $1$- and $\bar{1}$-labeled edges separate regions inside and outside of $D''$, e.g.,
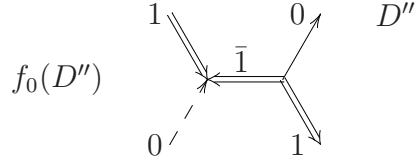
\[\begin{xy}
0;<1cm,0cm>:<.5cm,\halfrootthree cm>::
{\ar @{=>} (2,2);(3,1)}, {\ar @{-->} (3,0);(3,1)}, {\ar @{<=} (3,1);(4,1)}, {\ar @{->} (4,1);(4,2)}, {\ar @{=>} (4,1);(5,0)},
(1.8,2)*{1}, (3.7,2)*{0}, (3.3,1.3)*{\bar{1}}, (2.8,0)*{0}, (4.7,0)*{1},
(5,2)*{D''}, (1,1)*{f_0(D'')}
\end{xy}
\] The other cases are similar and the lemma follows.
\end{proof}

Let $D$ be a web, and let $e^*$ be the {\it {initial edge}} of $D$, namely the edge adjacent to the leftmost boundary vertex of $D$. Let $v^*$ be the other end of $e^*$. Construct the left cut $C^l := C^l_{e^*,v^*}$ and right cut $C^r:=C^r_{e^*,v^*}$, labeling their boundary endpoints $v^l$ and $v^r$, respectively. (We know these paths are disjoint after leaving vertex $v*$ and that they ultimately reach boundary vertices by Lemma \ref{lem:4}.)

We want to know the labelings of the edges on $C^l$ and $C^r$. Since $e^*$ is labeled with a $1$ and pointing toward $v^*$, we see that there are only two possibilities for the neighborhood of $v^*$. But by Lemma \ref{lem:1}, we know edge labelings are consistent with depth, and the edge to the left of $v^*$ separates two faces of depth $1$. Hence it is a $0$-labeled edge. Now by examination of the growth rules in Figure \ref{fig:pr3}, we see that any right turn from a downward-pointing $0$-edge takes us on an upward-pointing $1$-labeled edge. Any left turn from an upward-pointing $1$-edge leads to another downward-pointing $0$-edge and so on, as shown in Figure \ref{fig:pr8}. Because the path must have even length in order to end up on the boundary, we know that the final edge traversed is labeled with a $0$. Similarly, by examination of the local moves we have that $C^r$ alternates $\bar{1}0\bar{1}0\cdots$ upon leaving $v^*$, terminating at $v^r$, which, by parity considerations, must be labeled with $\bar{1}$.

We define $\Dl$ to be the collection of faces to the left of $C^l$ (when moving from $v^*$ to $v^l$). Similarly, $\Dr$ denotes the faces to the right of $C^r$ (notice that this includes the outer face $f_0$). Let $\Dm$ denote the faces to the right of $C^l$ and to the left of $C^r$. See Figure \ref{fig:pr8}.

\begin{figure}[h!]
\begin{center}
\[
\begin{xy}
0;<.75cm,0cm>:
(0,7); (16,7) **@{-},
{\ar @{->} (2.3,7);(2.3,6)}, {\ar @{->} (3.6,7);(3.6,6)}, {\ar @{->} (4.9,7);(4.9,6)}, {\ar @{->} (8.9,7);(8.9,6)}, {\ar @{<-} (8.9,6); (8.4,5.7)}, {\ar @{<-} (8.9,6); (9.4,5.7)}, {\ar @{->} (10.1,7);(10.1,6)}, {\ar @{<-} (10.1,6); (9.6,5.7)}, {\ar @{<-} (10.1,6); (10.6,5.7)},
{\ar @{->} (13.3,7);(13.3,6)}, {\ar @{->} (14.6,7);(14.6,6)},
{\ar@/_30pt/_1(1,7);(4.5,0)}, {\ar@/^10pt/^{\bar{1}}(10,0);(4.5,0)}, {\ar@/_5pt/_1(10,0);(13,0)},
(4.5,0) ="l1", (4.5,1.125)="l2", (5.25,1.875)="l3", (5.25,3)="l4", (6.25,4)="l5", (6.25,5.125)="l6", (7,5.875)="l7",
{\ar @{<-}_0 "l1"; "l2"}, {\ar @{->}_1 "l2"; "l3"}, {\ar @{<-}^0 "l3"; "l4"}, {\ar @{-->} "l4"; "l5"}, {\ar @{<-}_0 "l5"; "l6"}, {\ar @{->}_1 "l6"; "l7"}, {\ar @{<-}_0 "l7"; (7,7)}, {\ar @{->}_{\bar{1}} "l2";(3.75,2)}, {\ar @{<-}^{\bar{1}} "l3"; (6.25, 1.5)}, {\ar @{->}_{\bar{1}} "l4"; (4.5, 3.75)}, {\ar @{<-}_{\bar{1}} "l5"; (7.2, 3.5)}, {\ar @{->}_{\bar{1}} "l6"; (5.5, 5.5)},
{\ar @{<-}^{\bar{1}} "l7"; (8, 5.5)},
(10,0) ="r1", (10,1.125)="r2", (10.75,1.875)="r3", (10.75,3)="r4", (12.25,4)="r5", (11.25,5.125)="r6", (12,5.25)="r7",
{\ar @{->}_0 "r1"; "r2"}, {\ar @{<-}_{\bar{1}} "r2"; "r3"}, {\ar @{->}^0 "r3"; "r4"}, {\ar @{<--} "r4"; "r5"}, {\ar @{->}_0 "r5"; "r7"}, {\ar @{<-}_{\bar{1}} "r7"; (12,7)}, {\ar @{<-}_{1} "r2";(8.625,1.75)}, {\ar @{->}^{1} "r3"; (12.125, 1.5)}, {\ar @{<-}_{1} "r4"; (9.625, 3.625)}, {\ar @{->}_{1} "r5"; (13.25, 3.5)}, {\ar @{<-}_{1} "r7"; (10.625, 5.125)},
(1.5,6.75)="a1", (6.75,6.75)="a2", (4,.75)="a3",
{\ar @{--} "a1"; "a2"}, {\ar @{--} "a2"; "a3"}, {\ar @{--} @/^20pt/ "a3";"a1"},
(5,.75)="b1", (8,6.25)="b2", (8.9,5.5)="b3", (9.5,6.25)="b4", (10.1,5.5)="b5", (11.25,6.25)="b6", (9.25,.75)="b7",
{\ar @{--} "b1"; "b2"},  {\ar @{--} @/^3pt/ "b2"; "b3"}, {\ar @{--} @/^10pt/ "b3"; "b5"}, {\ar @{--} @/^3pt/ "b5"; "b6"}, {\ar @{--} @/^10pt/ "b6";"b7"}, {\ar @{--} @/^10pt/ "b7"; "b1"},
(1,7)*{\bullet}, (7,7)*{\bullet}, (12,7)*{\bullet}, (4.5,0)*{\bullet}, (2.3,7)*{\circ}, (3.6,7)*{\circ}, (4.9,7)*{\circ}, (8.9,7)*{\circ}, (10.1,7)*{\circ}, (13.3,7)*{\circ}, (14.6,7)*{\circ}, (2,-1)*{f_0}, (3,5)*{\Dl}, (8.5,4)*{\Dm}, (15,4)*{\Dr}, (3.5,1.25)*{L_1}, (7,0)*{M_1}, (11.5,.5)*{R_1}, (4.25,2.75)*{L_2}, (6.25,6.25)*{L_k}, (.5,5)*{e^*}, (4.5,-.5)*{v^*}, (7,7.5)*{v^l}, (12,7.5)*{v^r}
\end{xy}
\]
\end{center}
\caption{}\label{fig:pr8}
\end{figure}

\begin{lemma} \label{lem:5}
Let $D$ be a web. After moving the leftmost boundary vertex to the right,
\begin{enumerate}
\item the depth of every face in $\Dl$ decreases by 1,
\item the depth of every face in $\Dr$ increases by 1, and
\item the depth of every face in $\Dm$ remains unchanged.
\end{enumerate}
\end{lemma}

\begin{proof}
Let $L_1$ denote the face separated from the outer face by $e^*$. This face will be the outer face once the leftmost boundary vertex moves to the right. Let $L_2, \ldots, L_k$ denote the other faces of $\Dl$ that border the left cut. By examining the edge labels (which by Lemma \ref{lem:1} are consistent with depth) every face $L_i$ has a minimal path to $f_0$ that passes through $L_1$. Thus, any face in $\Dl$ has a minimal path to $f_0$ that goes through $L_1$. Claim (1) then follows.

By examining the faces on the boundary of $\Dm$, we see that no face in $\Dm$ has a minimal length path through $L_1$, but they all have such a path through $M_1$. Since $M_1$ is a neighbor to both $f_0$ and $L_1$, this implies (3). A similar argument shows that any face in $\Dr$ is closer to $R_1$ than to $M_1$, and (2) follows.
\end{proof}

According to Lemma \ref{lem:5}, we have now established that all the vertices to the left of $v^l$ keep their labels when we move the leftmost vertex to the right, while $v^l$'s label changes to a $1$. Likewise all the labels between $v^l$ and $v^r$ are the same, but $v^r$ has changed to a $0$. If $a$ is the position of $v^l$ in the word $w = w(D)$, and $b$ is the position of $v^r$, then we have $w' = w(p(D))$ given by: \[ w' = w_2 \cdots w_{a-1} 1 w_{a+1} \cdots w_{b-1} 0 w_{b+1} \cdots w_{3n} \bar{1}.\]
All that remains is to verify that $a$ is the first position where letters $1$ and $0$ balance in $w$, and $b$ is the first position where letters $0$ and $\bar{1}$ balance. This is the content of Lemma \ref{lem:6}.

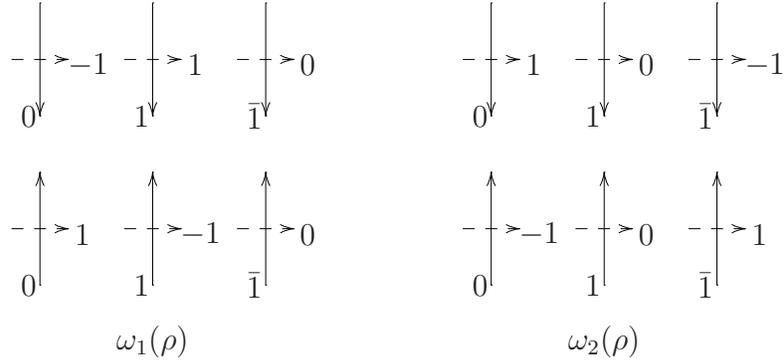
\begin{figure}[h!]
\begin{center}
\[
\begin{xy}
0;<.75cm,0cm>:
{\ar @{->} (0,1);(0,3)}, {\ar @{-->} (-.5,2);(.5,2)},
(-.2,1)*{0}, (.75,1.9)*{1},
{\ar @{->} (2,1);(2,3)}, {\ar @{-->} (1.5,2);(2.5,2)},
(1.8,1)*{1}, (2.85,1.9)*{-1},
{\ar @{->} (4,1);(4,3)}, {\ar @{-->} (3.5,2);(4.5,2)},
(3.8,1)*{\bar{1}}, (4.75,1.9)*{0},
{\ar @{<-} (0,4);(0,6)}, {\ar @{-->} (-.5,5);(.5,5)},
(-.2,4)*{0}, (.85,4.9)*{-1},
{\ar @{<-} (2,4);(2,6)}, {\ar @{-->} (1.5,5);(2.5,5)},
(1.8,4)*{1}, (2.75,4.9)*{1},
{\ar @{<-} (4,4);(4,6)}, {\ar @{-->} (3.5,5);(4.5,5)},
(3.8,4)*{\bar{1}}, (4.75,4.9)*{0},
{\ar @{->} (8,1);(8,3)}, {\ar @{-->} (7.5,2);(8.5,2)},
(7.8,1)*{0}, (8.85,1.9)*{-1},
{\ar @{->} (10,1);(10,3)}, {\ar @{-->} (9.5,2);(10.5,2)},
(9.8,1)*{1}, (10.75,1.9)*{0},
{\ar @{->} (12,1);(12,3)}, {\ar @{-->} (11.5,2);(12.5,2)},
(11.8,1)*{\bar{1}}, (12.75,1.9)*{1},
{\ar @{<-} (8,4);(8,6)}, {\ar @{-->} (7.5,5);(8.5,5)},
(7.8,4)*{0}, (8.75,4.9)*{1},
{\ar @{<-} (10,4);(10,6)}, {\ar @{-->} (9.5,5);(10.5,5)},
(9.8,4)*{1}, (10.75,4.9)*{0},
{\ar @{<-} (12,4);(12,6)}, {\ar @{-->} (11.5,5);(12.5,5)},
(11.8,4)*{\bar{1}}, (12.85,4.9)*{-1},
(2,0)*{\omega_1(\rho)},
(10,0)*{\omega_2(\rho)}
\end{xy}
\]
\end{center}
\caption{Definition of weights $\omega_1$ and $\omega_2$.}\label{fig:pr4}
\end{figure}

But first, we need one more tool. Let a (directed) curved line $\ell$ intersect a web $D$ so that it does not pass through any vertices. To each point $\rho \in \ell$ that intersects an edge of $D$ we assign weights $\omega_1(\rho)$ and $\omega_2(\rho)$ according to the rules shown in Figure \ref{fig:pr4}.
Here the dashed line denotes $\ell$ and the numbers next to it denote the values of $\omega_1(\rho)$ and $\omega_2(\rho)$, respectively. Finally, we let $\omega_i(\ell) = \sum_\rho \omega_i(\rho)$ where the sum is taken over all intersections of $\ell$ with $D$.

The following lemma shows that for fixed starting and ending points, these statistics are independent of the path chosen.

\begin{lemma} \label{lem:2}
The values of $\omega_1(\ell)$ and $\omega_2(\ell)$ depend only on the endpoints of $\ell$ and not on the exact path it takes. 
\end{lemma}

\begin{proof}
 The statement follows from verification of the local moves as shown on Figure \ref{fig:pr6}, where $\{i,j,k\} = \{1,0,\bar 1\}$. Note that according to growth rules edges adjacent to any internal vertex are labelled this way. For example, take $i=1, j=0, k=\bar{1}$, the arrows oriented towards the central vertex and the dashed line $l$ directed eastward. Then by the rules on Figure \ref{fig:pr4} we see $\omega_1(\ell) = 1$ and $\omega_2(\ell) = 0$, regardless of the path we take. All other cases are similar.
 
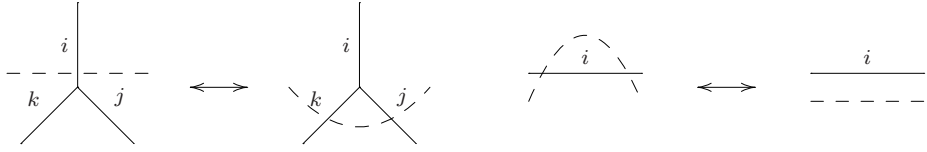
\begin{figure}[h!]
\begin{center}
\[
\begin{xy}
0;<.75cm,0cm>:
{\ar @{-}_k (1,1); (0,0)}, {\ar @{-}^j (1,1); (2,0)}, {\ar @{-}^i (1,1); (1,2.5)}, {\ar @{--} (-.25,1.25); (2.25,1.25)},
{\ar @{<->} (3,1);(4,1)}, {\ar @{-}_k (6,1); (5,0)}, {\ar @{-}^j (6,1); (7,0)}, {\ar @{-}^i (6,1); (6,2.5)}, {\ar @{--} @/_15pt/ (4.75,1); (7.25,1)},
{\ar @{-}^i (9,1.25); (11,1.25)}, {\ar @{--} @/^25pt/ (9,.75); (11,.75)}, {\ar @{<->} (12,1); (13,1)}, {\ar @{-}^i (14,1.25); (16,1.25)}, {\ar @{--} (14,.75); (16,.75)} 
\end{xy}
\]
\end{center}
\caption{Weight-preserving moves.}\label{fig:pr6}
\end{figure}

Clearly any route $\ell$ might take between any two fixed endpoints can be transformed into any other route by a sequence of such moves. This proves the lemma.
\end{proof}

Now we are ready to show that $v^l$ and $v^r$ are located at the proper positions in the word $w(D)$, establishing Theorem \ref{thm:prom}.

\begin{lemma} \label{lem:6}
Let $v^l$ and $v^r$ be vertices in $D$ defined as before.
\begin{enumerate}
\item Among the labels of vertices preceding $v^l$ (inclusively) there are as many letters $1$ as $0$, i.e., any path $\ell$ from $f_0$ to the face to the left of $v^l$ has $\omega_1(\ell) = 1$. Further, $v^l$ is the leftmost vertex with this property.
\item Among the labels of vertices preceding $v^r$ (inclusively) there are as many letters $0$ as $\bar{1}$, i.e., any path $\ell$ from $f_0$ to the face to the left of $v^r$ has $\omega_2(\ell) = 1$. Further, $v^r$ is the leftmost vertex to the right of $v^l$ with this property.
\end{enumerate}
\end{lemma}

\begin{proof}
Using Lemma \ref{lem:2}, we see that the $\omega_1(\ell) = 1$ for any path $\ell$ from $f_0$ to the boundary face to the left of the $v^l$ ($L_k$ in Figure \ref{fig:pr8}). This value is easily computed by taking a path just to the left of $C^l$. 

Now for any face $f$ in $\Dl$, define $\omega_1(f)$ as the value of $\omega_1(\ell)$ for any path $\ell$ from $f_0$ to $f$. Similarly, let $\omega'_1(f)$ denote the weight $\omega_1(\ell')$ of any path $\ell'$ from $L_1$ to $f$. Then we have $\omega'_1(f) = \omega_1(f) -1$.

To see that $v^l$ is the leftmost vertex with the desired property, consider the web $D_L$ formed by taking $\Dl - L_1$ (enclosed by a dashed line in Figure \ref{fig:pr8}). This is a web with all inward pointing edges, and so by Theorem \ref{thm:kkh} it must be dominant, i.e., the word $w(D_L)$ is Yamanouchi. In particular, $\omega'_1(f) \geq 0$ for any face $f$ on the boundary of $L'$.

If there was a position to the left of $v^l$ where letters $1$ and $0$ balance in $w(D)$, then there would be a face $f$ on the boundary of both $D$ and $D_L$ where $\omega_1(f) = 0$. But then $\omega'_1(f) = -1$, a contradiction.

For $v^r$, the reasoning is similar. Walking along the left side of the right cut allows us to compute $\omega_2(\ell) = 1$ for any path from $f_0$ to the boundary face to the left of $v^r$.

Now for a face $f$ in $\Dm$, we let $\omega_2(f)$ denote the value of $\omega_2(\ell)$ for any path $\ell$ from $f_0$ to $f$; $\omega'_2(f)$ denotes the weight $\omega_2(\ell')$ of any path $\ell'$ from $M_1$ to $f$. Clearly, $\omega'_2(f) = \omega_2(f)-1$.

Consider the web $D_M$ formed by starting in $M_1$, cutting along the right side of $C^l$, zig-zagging along the boundary of $D$ to collect edges of the same parity, then passing down the left side of $C^r$ (outlined in Figure \ref{fig:pr8}). To be more precise, along the boundary of $D$ we ``zig-zag" in two ways. In type I, two consecutive edges join up, in which case we take the third edge connected to them. In type II, one edge has two non-boundary branches that we pass through. Along the left and right cuts, all moves are type I.
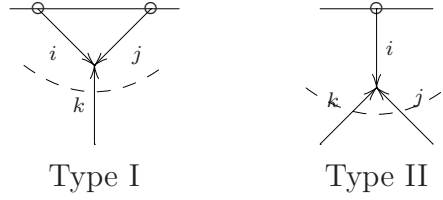
\[
\begin{xy}
0;<.75cm,0cm>:
(-.5,2); (2.5,2) **@{-},
(0,2)="a", (1,1)="b", (2,2)="c", (1,-.4)="d",
"a"*{\circ}, "c"*{\circ}, {\ar @{->}_i "a";"b"}, {\ar @{->}^j "c";"b"}, {\ar @{->}^k "d";"b"}, {\ar @{--} @/_10pt/ (-.25,1); (2.25,1)},
(1,-1)*{\mbox{Type I}},
(5,2); (7,2) **@{-},
(6,2)="e", (6,.6)="f", (5,-.4)="g", (7,-.4)="h",
"e"*{\circ}, {\ar @{->}^i "e"; "f"}, {\ar @{->}^k "g"; "f"}, {\ar @{->}_j "h"; "f"}, {\ar @{--} @/_10pt/ (4.75,.6); (7.25,.6)},
(6,-1)*{\mbox{Type II}}
\end{xy}
\]

Since all the boundary edges have the same orientation, we can apply Theorem \ref{thm:kkh} to conclude that the word $w(D_M)$ is Yamanouchi; in particular, $\omega'_2(f) \geq 0$ for any boundary face $f$ of $M'$.

If there is a boundary vertex strictly between $v^l$ and $v^r$ where letters $0$ and $\bar{1}$ balance, this means there is a face $f$ on the boundary of both $D$ and $D_M$ where $\omega_2(f) = 0$, but then this implies that $\omega'_2(f) = -1$, a contradiction. (Such a face cannot occur in $D - D_M$ since we cannot have $i = \bar{1}$ in a type I crossing.)

This completes the proof of the lemma.
\end{proof}

We have now proved Theorem \ref{thm:prom} as well.

\section{Application to cyclic sieving}\label{sec:sieving}

Let $W_{n}^{(2)}$ denote the $\mathbb{C}$-vector space of irreducible $A_1$ webs with $2n$ boundary vertices, that is, noncrossing matchings on $[2n]$.  Similarly, let $W_{n}^{(3)}$ denote the $\mathbb{C}$-vector space spanned by the set $\mathfrak{M}_n$ of irreducible $A_2$ webs with $3n$ $``-"$s on the boundary.  

We define actions of $S_{2n}$ on $W_{n}^{(2)}$ and $S_{3n}$ on $W_{n}^{(3)}$ as follows.  For an $A_1$ web $E$ and an index $i \in [2n-1]$, define $t^{(2)}_i \cdot E$ to be the element of $W_{n}^{(2)}$ obtained by attaching an uncrossing ''\noXv'' at vertices $i$ and $i+1$ to the diagram of $E$. Here we apply the $A_1$ spider reduction rule if nesessery, that is we replace any resulting loop with a factor of $-2$.  Define the action of the Coxeter generator $s_i = (i,i+1) \in S_{2n}$ on $W_{n}^{(2)}$ by $s_i \cdot E := E + t^{(2)}_i \cdot E$ for all noncrossing matchings $E$, extended linearly.

\begin{figure}[h!]
\begin{center}
\input{pr9.pstex_t}
\end{center}
\caption{}\label{fig:pr9}
\end{figure}

Similarly, for any irredicible $A_2$ web $D \in \mathfrak{M}_n$ and an index $i \in [3n-1]$, set $t^{(3)}_i \cdot D$ equal to the element of $W_{n}^{(3)}$ obtained by attaching an uncrossing ``\noXvsp`` at indices $i$ and $i+1$.  Proceed expressing the resulting $A_2$ web $t^{(3)}_i \cdot D$ as a linear combination of irreducible $A_2$ webs via the spider reduction rules.  We define the action of the Coxeter generator $s_i = (i,i+1)$ on $W_{n}^{(3)}$ by $s_i \cdot D := D + t_i^{(3)}D$.  

\begin{figure}[h!]
\begin{center}
\input{pr10.pstex_t}
\end{center}
\caption{}\label{fig:pr10}
\end{figure}

\begin{lemma}
The actions of the Coxeter generators on $W_{n}^{(2)}$ and $W_{n}^{(3)}$ defined above extend to actions of the appropriate symmetric groups to make $W_{n}^{(2)}$ and $W_{n}^{(3)}$ modules over $S_{2n}$ and $S_{3n}$, respectively. 
\end{lemma}

\begin{proof}
We must verify that the Coxeter relations are satisfied.  This is an easy exercise involving the spider reduction rules in the case of $W_n^{(3)}$ and an easier exercise involving the relation that a closed loop yields a factor of $-2$ in the case of $W_n^{(2)}$.
\end{proof} 

In fact, the resulting action is the action of certain quotients of the group algebra of the symmetric group. Namely, the {\it {Temperley-Lieb algebra}} for $A_1$, and a {\it {Temperley-Lieb-Martin algebra}} \cite{M} for $A_2$, cf. \cite{P}.

Next, we identify $W_{n}^{(2)}$ and $W_{n}^{(3)}$ as irreducible modules over $S_{2n}$ and $S_{3n}$, respectively.

\begin{lemma}
\begin{enumerate}
 \item $W_{n}^{(2)}$ is an  irreducible $S_{2n}$-module of shape $(n,n)$.\\
 \item $W_{n}^{(3)}$ is an irreducible $S_{3n}$-module of shape $(n,n,n)$.
\end{enumerate}
\end{lemma}

\begin{proof}
Let $\rho^{(2)} : \mathbb{C}[S_{2n}] \rightarrow End(W_{n}^{(2)})$ and 
$\rho^{(3)} : \mathbb{C}[S_{3n}] \rightarrow End(W_{n}^{(3)})$ denote the algebra homomorphisms which define the module structure for $W_{n}^{(2)}$ and $W_{n}^{(3)}$.  For any subset $X \subseteq S_n$, define $[X]_-$ to be the group algebra element given by 
$$[X]_- = \sum_{x \in X} sgn(x) x.$$  For any partition $\lambda \vdash n$, define $S_{\lambda}$ to be the Young subgroup of $S_n$ indexed by $\lambda$.  That is, $S_{\lambda}$ is the subgroup of $S_n$ which fixes setwise the sets $\{1,2,\dots, \lambda_1 \}$, $\{\lambda_1+1,\lambda_1+2,\dots,\lambda_1+\lambda_2 \}, \dots$.   

Since the action of the symmetric group on $W_n^{(2)}$ and $W_n^{(3)}$ factors through the Temperley-Lieb algebra and the Temperley-Lieb-Martin algebra, one concludes that the irreducible components cannot have more than $2$ and $3$ rows correspondingly. On the other hand, it is easy to show that $[S_{(n^2)}]_-$ and $[S_{(n^3)}]_-$ do indeed act nontrivially on $W_n^{(2)}$ and $W_n^{(3)}$.  Since $$[S_{\lambda}]_- {\mathbb C} S_n = Ind \uparrow_{S_{\lambda}}^{S_n}(\bf 1')$$ as a left ${\mathbb C} S_n$-module, where $\bf 1'$ is the alternating representation, we can use the fact that the Kostka matrix is upper triangular with respect to dominance order. We conclude that $W_{n}^{(2)}$ and $W_{n}^{(3)}$ contain irreducible components smaller than or equal to the corresponding rectangular shapes. The only shape that is not larger than a $k$ by $n$ rectangle in dominance order but has at most $k$ rows is the rectangle itself. Finally, dimension count shows that this irreducible occurs in $W_{n}^{(k)}$ exactly once for $k=2,3$ correspondingly, while others do not occur.
\end{proof}

We are almost ready to give a proof of the desired CSP, but first we want to have a more compact way of realizing the action of the Coxeter generators on webs. 

To do so, we extend the notion of webs to allow crossings as follows. In $A_1$, a crossing should be understood as the state sum
\begin{center}
\begin{picture}(0,0)%
\includegraphics{pr11.pstex}%
\end{picture}%
\setlength{\unitlength}{2486sp}%
\begingroup\makeatletter\ifx\SetFigFontNFSS\undefined%
\gdef\SetFigFontNFSS#1#2#3#4#5{%
  \reset@font\fontsize{#1}{#2pt}%
  \fontfamily{#3}\fontseries{#4}\fontshape{#5}%
  \selectfont}%
\fi\endgroup%
\begin{picture}(4074,699)(1564,-523)
\put(2701,-331){\makebox(0,0)[lb]{\smash{{\SetFigFontNFSS{11}{13.2}{\familydefault}{\mddefault}{\updefault}{\color[rgb]{0,0,0}$=$}%
}}}}
\put(4276,-331){\makebox(0,0)[lb]{\smash{{\SetFigFontNFSS{11}{13.2}{\familydefault}{\mddefault}{\updefault}{\color[rgb]{0,0,0}$+$}%
}}}}
\end{picture}%

\end{center}
and for $A_2$ webs as the state sum 
\begin{center}
\begin{picture}(0,0)%
\includegraphics{pr12.pstex}%
\end{picture}%
\setlength{\unitlength}{2486sp}%
\begingroup\makeatletter\ifx\SetFigFontNFSS\undefined%
\gdef\SetFigFontNFSS#1#2#3#4#5{%
  \reset@font\fontsize{#1}{#2pt}%
  \fontfamily{#3}\fontseries{#4}\fontshape{#5}%
  \selectfont}%
\fi\endgroup%
\begin{picture}(4074,759)(1564,-553)
\put(2701,-331){\makebox(0,0)[lb]{\smash{{\SetFigFontNFSS{11}{13.2}{\familydefault}{\mddefault}{\updefault}{\color[rgb]{0,0,0}$=$}%
}}}}
\put(4276,-331){\makebox(0,0)[lb]{\smash{{\SetFigFontNFSS{11}{13.2}{\familydefault}{\mddefault}{\updefault}{\color[rgb]{0,0,0}$+$}%
}}}}
\end{picture}%

\end{center}
Now we see that a simple transposition $s_i$ simply introduces a crossing between boundary vertices $i$ and $i+1$.

With this viewpoint, it is straightforward to check that the following Reidemester-type moves can be performed for $A_1$ and $A_2$ webs correspondingly. In the $A_2$ case one should interpret unoriented edges as either of the two possible orientations.
\begin{center}
\begin{picture}(0,0)%
\includegraphics{pr13.pstex}%
\end{picture}%
\setlength{\unitlength}{1657sp}%
\begingroup\makeatletter\ifx\SetFigFontNFSS\undefined%
\gdef\SetFigFontNFSS#1#2#3#4#5{%
  \reset@font\fontsize{#1}{#2pt}%
  \fontfamily{#3}\fontseries{#4}\fontshape{#5}%
  \selectfont}%
\fi\endgroup%
\begin{picture}(10595,4070)(1343,-5244)
\put(5716,-1996){\makebox(0,0)[lb]{\smash{{\SetFigFontNFSS{7}{8.4}{\familydefault}{\mddefault}{\updefault}{\color[rgb]{0,0,0}$=-$}%
}}}}
\put(2206,-1996){\makebox(0,0)[lb]{\smash{{\SetFigFontNFSS{7}{8.4}{\familydefault}{\mddefault}{\updefault}{\color[rgb]{0,0,0}$=$}%
}}}}
\put(2206,-4651){\makebox(0,0)[lb]{\smash{{\SetFigFontNFSS{7}{8.4}{\familydefault}{\mddefault}{\updefault}{\color[rgb]{0,0,0}$=$}%
}}}}
\put(5761,-4606){\makebox(0,0)[lb]{\smash{{\SetFigFontNFSS{7}{8.4}{\familydefault}{\mddefault}{\updefault}{\color[rgb]{0,0,0}$=$}%
}}}}
\put(10036,-4651){\makebox(0,0)[lb]{\smash{{\SetFigFontNFSS{7}{8.4}{\familydefault}{\mddefault}{\updefault}{\color[rgb]{0,0,0}$=-$}%
}}}}
\end{picture}%

\end{center}

Let $N = bn$, $b =2,3$. We will now relate the action of web rotation to the action of the long cycle  $c = (1 2 \cdots N)$ in $S_N$.

\begin{lemma}\label{lem:rotcyc}
For $b = 2,3$, the action of rotation of an $A_{b-1}$ web $D$ is, up to sign, the action of the long cycle, i.e., \[ p(D) = (-1)^{b-1} c\cdot D.\]
\end{lemma}

\begin{proof}
By iterating the crossings corresponding to the Coxeter generators, we see the long cycle $c = (1 2 \cdots N) = s_{N-1} \cdots s_2 s_1$ acts as the ``whirl" shown below. 
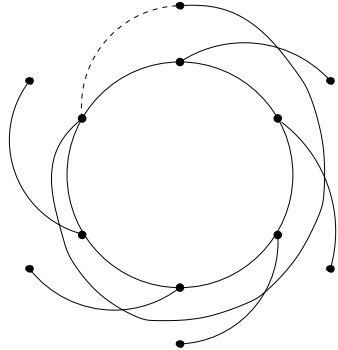
\begin{center}
\begin{picture}(0,0)%
\includegraphics{pr14.pstex}%
\end{picture}%
\setlength{\unitlength}{2072sp}%
\begingroup\makeatletter\ifx\SetFigFontNFSS\undefined%
\gdef\SetFigFontNFSS#1#2#3#4#5{%
  \reset@font\fontsize{#1}{#2pt}%
  \fontfamily{#3}\fontseries{#4}\fontshape{#5}%
  \selectfont}%
\fi\endgroup%
\begin{picture}(3922,4156)(3126,-5514)
\end{picture}%

\end{center}
We want to use Reidemeister moves to pull the long string that wraps around into the position shown by dashed line, forming the diagram for rotation. This can clearly be done. Furthermore, in the $A_1$ case one needs to apply the sign changing transformation exactly once, while in $A_2$ one needs to apply the sign-changing transformation exactly twice. (The endpoints of the strings are fixed.)
\end{proof}

\begin{proposition}
Let $\lambda \vdash N = bn$ be a rectangle with $b=2$ or $3$ rows and let $C = \mathbb{Z} / N \mathbb{Z}$ act on $X = SYT(\lambda)$ by promotion.  Then the triple $(X, C, X(q))$ exhibits the cyclic sieving phenomenon, where $X(q)$ is as in Theorem \ref{thm:cs}.
\end{proposition} 

\begin{proof} 
By Lemma \ref{lem:rotcyc}, for any $d \geq 0$, the number of webs on $N$ vertices fixed by $d$ rotations is equal to the value of the irreducible character $\chi^{(n,n)}$  or $\chi^{(n,n,n)}$ of $S_N$ evaluated on the permutation $c_N^d$.  In order to get our cyclic sieving result, we need to relate this character evaluation to a polynomial evaluation.  To do this, we use Springer's theory of regular elements \cite{Sp}.  For $W$ a finite complex reflection group, an element $w \in W$ is called \emph{regular} if there exists an eigenvector $v$ for $w$ in the reflection representation of $W$ such that $v$ does not lie on any of the reflecting hyperplanes for the reflections in $W$.  

Keeping the notation of the previous paragraph, let $\chi^{\lambda}$ be an irreducible character of $W$.  We can associate to $\chi^{\lambda}$ a polynomial called the \emph{fake degree polynomial} as follows.  Letting $V$ denote the reflection representation of $W$, let $\mathbb{C}[V]$ denote the ring of polynomial valued functions on $V$ and let $\mathbb{C}[V]^W_+$ denote the subring of those functions which are invariant under the action of $W$.  The quotient $\mathbb{C}[V] / \mathbb{C}[V]^W_+$ carries an action of $W$ which is graded.  Define the fake degree polynomial $f^{\lambda}(q) = \sum_{i \geq 0} a_i q^i$ by letting $a_i$ be the multiplicity of $\chi^{\lambda}$ in the $i$-th graded piece of this representation.
Springer showed that if $w$ is a regular element of $W$ and $v$ is an associated eigenvector and $w \cdot v = \omega v$, we have that $\chi^{\lambda}(w) = f^{\lambda}(\omega)$.    

We apply Springer's result to the case of $W = S_N$ to get our desired cyclic sieving phenomenon.
It is easy to see that $c_N^d$ is a regular element of $S_N$ for all $d$.  Moreover, it is possible to show that for any partition $\lambda \vdash N$ the fake degree polynomial for the irreducible representation of $S_N$ with shape $\lambda$ has the following form:
\begin{equation*}
f^{\lambda}(q) = q^{-\kappa(\lambda)} \frac{[N]_q!}{\Pi_{(i,j) \in \lambda} [h_{ij}]_q},
\end{equation*}
where $\kappa(\lambda) = 0 \lambda_1 + 1 \lambda_2 + 2 \lambda_3 + \cdots$.  
For the $A_1$ case, assume $\lambda = (n,n)$ has two rows.  Then, $\kappa(\lambda) = n$ and if $\zeta$ is a primitive $N^{th}$ root of unity with $2n = N$, then $(\zeta)^{-nd} = (-1)^d$ for all $d \geq 0$.  On the other hand, by Lemma \ref{lem:rotcyc} we also have that $\chi^{(n,n)}(c_N^d)$ is equal to $(-1)^d$ times the number of elements of $SYT((n,n))$ fixed under $d$ iterations of promotion.  The desired CSP follows.  For the $A_2$ case, notice that if $\lambda = (n,n,n)$ has three rows and $3n = N$ and $\zeta$ is a primitive $N^{th}$ root of unity, we have that $\kappa(\lambda) = 3n$ and 
$(\zeta)^{-3nd} = 1$ for all $d \geq 0$.  On the other hand, in this case $\chi^{\lambda}(c_N^d)$ is equal to the number of elements of $SYT((n,n,n))$ fixed by $d$ iterations of promotion, completing the proof.
\end{proof}

\subsection{Enumeration of web orbits}

We can now extract the number of $A_2$ webs fixed by any given number of rotations, $d|3n$, by taking $q \to e^{2\pi i/d}$ in \[f^{(n,n,n)}(q) = \frac{[3n]_q! [2]_q}{[n]_q![n+1]_q![n+2]_q!}.\] These numbers are something that we have no way to compute other than via the CSP, though formula \eqref{eq:rot} suggests that a more direct argument may exist.

\begin{proposition}\label{prp:numrot}
For $n\geq 3$, the number of webs fixed by $3n/d$ rotations is the multinomial coefficient

\begin{equation}\label{eq:rot}
\binom{ 3n/d }{ \lfloor n/d \rfloor, \lfloor (n+1)/d \rfloor, \lfloor (n+2)/d\rfloor} = \frac{(3n/d)!}{ \lfloor n/d \rfloor! \lfloor (n+1)/d \rfloor! \lfloor (n+2)/d\rfloor!},
\end{equation}
if $d=3$ or $d|n$, zero otherwise.

\end{proposition}

\begin{remark} The similar exercise for noncrossing matchings is, under a bijection with triangulations of polygons, handled in \cite[Theorem 7.1]{RSW}.
\end{remark}

\begin{remark}
The condition that $d=3$ or $d|n$ means that many proper divisors of $3n$ will not fix any webs. For instance, there are no webs with 24 vertices fixed by four rotations since $d=6$ does not divide $n=8$.
\end{remark}

\begin{remark}
With $d =2$, $n = 2k$, equation \eqref{eq:rot} gives $\frac{(3k)!}{k!k!(k+1)!}$, which is not strictly speaking a multinomial coefficient since $k + k + (k+1) \neq 3k$. For $d > 2$, the number is a true multinomial.
\end{remark}

\begin{corollary}
For $n\geq 3$, there are six webs on $3n$ vertices fixed by three rotations; those in the orbit of $w = (123)^n$, and those in the orbit of $w'=11122(132)^{n-3}2333$.
\end{corollary}

\begin{proof}
Taking $d = n$ in \eqref{eq:rot} we get $\binom{3}{1,1,1} = 6$ webs fixed by $3n/d = 3$ rotations. By considering how promotion acts on Yamanouchi words, it is not difficult to verify that the webs $w$ and $w'$ have the desired orbits of size three.
\end{proof}

\begin{proof}[Proof of Proposition \ref{prp:numrot}]
We proceed by evaluation of the hook length formula at appropriate roots of unity, i.e., primitive $d$th roots of unity, where $d|3n$.

Let $\zeta = e^{2\pi i/d}$. We now apply the following rules (throughout this proof we abbreviate $[m]_q$ by $[m]$): \[\lim_{q \to \zeta} \frac{[m_1]}{[m_2]} = \begin{cases}
\frac{m_1}{m_2} &\mbox{if } m_1 \equiv m_2 \equiv 0 \mod d,\\
1 &\mbox{if } m_1 \equiv m_2 \not \equiv 0 \mod d,
\end{cases}
\]
and \[ \lim_{q \to \zeta} [m] = 0 \mbox{ if and only if } d|m.\]

For any $d|3n$ we have 
\begin{align*}
\lim_{q \to \zeta} f^{(n,n,n)}(q) &=\frac{ [3n]![2]}{[n]![n+1]![n+2]!} \\
 &= \lim_{q \to \zeta} \frac{ [3n] \cdots [n+3]}{[n] \cdots [2][n+1][n]\cdots [3]}\frac{[n+2][n+1] \cdots [2][2]}{[2][n+2] \cdots [2]}\\
&= \lim_{q \to \zeta} \frac{ [3n] \cdots [n+3]}{[n] \cdots [2][n+1][n]\cdots [3]}.
\end{align*}

If $3<d|3n$ but $d$ does not divide $n$, then there are always more terms $[m]$ in the numerator for which $d|m$ than in the denominator, forcing $\lim_{q \to \zeta} f^{(n,n,n)}(q) = 0$ in this case.

For $d|n=dk$, $3\leq d \leq n$, we have
\begin{align*}
\lim_{q \to \zeta} f^{(n,n,n)}(q) &= \lim_{q \to \zeta} \frac{ [3n] }{[n]}\cdots \frac{[2n+2]}{[2]}  \frac{[2n+1]}{[n+1]}\cdots \frac{[n+3]}{[3]} \\
&= \frac{3n}{n}\frac{(3n-d)}{(n-d)} \cdots \frac{(2n+d)}{d}\frac{2n}{n}\frac{(2n-d)}{(n-d)}\cdots \frac{(n+d)}{d}\\
&= \frac{ d( (3k)(3k-1) \cdots (2k+1)(2k)(2k-1) \cdots (k+1) )}{d( k(k-1)\cdots 1 k(k-1) \cdots 1)} =\frac{(3k)!}{k!k!k!}.
\end{align*}

The case for $d=2|n$ is similar, as are the cases when $d=3$ does not divide $n$.
\end{proof}

\begin{example}
For $n=4$, with $\zeta = e^{\pi i/6}$, we have 
\[
\begin{array}{c c c c}
f^{(4,4,4)}(1) = 462, & f^{(4,4,4)}(\zeta) = 0, & f^{(4,4,4)}(\zeta^2) = 0, & f^{(4,4,4)}(\zeta^3) = 6, \\
f^{(4,4,4)}(\zeta^4) = 12, & f^{(4,4,4)}(\zeta^5) = 0, & f^{(4,4,4)}(\zeta^6) = 30, & f^{(4,4,4)}(\zeta^7) = 0, \\
f^{(4,4,4)}(\zeta^8) = 12, & f^{(4,4,4)}(\zeta^9) = 6, & f^{(4,4,4)}(\zeta^{10}) = 0, & f^{(4,4,4)}(\zeta^{11}) = 0. 
\end{array}
\]
Further, we can deduce the sizes of the orbits of promotion/rotation. Let $o_k$ denote the number of $k$-orbits. Since six webs are fixed by three rotations, $o_3 = 2$. Since twelve webs are fixed by four rotations (and none are fixed by two rotations), we get $o_4 = 3$. We get $o_6 = 4$ because six of the thirty webs fixed by six rotations live in $3$-orbits. Proceeding, we see that the remaining $462 - 42 = 420$ webs must live in $12$-orbits, and $o_{12} = 35$.

Such an enumeration of $k$-orbits is possible whenever the CSP is present. See \cite{RSW}.
\end{example}

\section{Concluding remarks}\label{sec:concl}

There are several potential avenues for further study of the questions raised in this paper. Perhaps the most obvious of these is whether our approach can be used to prove Theorems \ref{thm:ordern} and \ref{thm:cs} for arbitrary rectangles. Unfortunately, this seems to be rather unlikely for the moment. The theory of spiders has yet to be generalized beyond the rank 2 case, though Kim \cite{Ki} has conjectured relations for an $A_3$ spider, and Jeong and Kim \cite[Theorem 2.4]{JK} show that $A_n$ ``webs" can be defined as planar graphs. Generalizing our approach in this way will depend on a concrete description for $A_n$ webs.

Another idea for generalization is to examine the other spiders that are well-understood; namely the $B_2$ and $G_2$ spiders found in \cite{K}, as well as $B_3$ case considered in \cite{We}. For these we can ask two questions: does their rotation correspond to some known generalization of promotion? and, do these webs admit a cyclic sieving phenomenon? Potentially related is Haiman's theorem for generalized staircases \cite[Theorem 4.4]{H}, which classifies the shifted shapes $\lambda$ for which $|\lambda|$ promotions fixes all tableaux. However, it may be that relation of promotion and the cyclic rotation of webs found here is a type $A$ phenomenon. In this case one can still look for cyclic sieving phenomena for webs of other types.

\end{document}